\documentclass[11pt]{article}
\usepackage{amsmath,amsfonts,amssymb}
\usepackage{bbm}
\usepackage{cases}
\usepackage[hidelinks]{hyperref}

\def\|{\vert\vert}

\def\<{{\langle}}
\def\>{{\rangle}}

\newtheorem{theorem}{Theorem}[section]
\newtheorem{lemma}[theorem]{Lemma}

\newtheorem{remark}[theorem]{Remark}

\newtheorem{definition}[theorem]{Definition}

\begin{document}

\title{\vspace*{-1.5cm}
Stochastic heat equations with logarithmic nonlinearity}

\author{Shijie Shang$^{1}$ and
Tusheng Zhang$^{2}$}
\footnotetext[1]{\, School of Mathematics, University of Science and Technology of China, Hefei, China.}
\footnotetext[2]{\, School of Mathematics, University of Manchester, Oxford Road, Manchester
M13 9PL, England, U.K. Email: tusheng.zhang@manchester.ac.uk}
\maketitle

\begin{abstract}
In this paper, we establish the existence and uniqueness of solutions to stochastic heat equations with logarithmic nonlinearity driven by Brownian motion on a bounded domain $D$ in the setting of $L^2(D)$ space. The result is valid for all initial values in $L^2(D)$. The logarithmic Sobolev inequality plays an important role.
\end{abstract}


\medskip


\section{Introduction}

In this paper, we study the stochastic heat equation with a logarithmic nonlinear term driven by  Brownian motion, which is given as follows:
\begin{numcases}{}
 du(t,x) = \Delta u(t,x) dt+ u(t,x)\log|u(t,x)| dt+ \sigma(u(t,x))dB_t, \   t>0, x\in D , \nonumber\\
 u(t,x)=0, \quad t>0, x\in\partial D , \nonumber\\
\label{1.a} u(0,x)=u_0(x), \quad  x\in D ,
\end{numcases}
where $D$ is a bounded domain of $\mathbb{R}^d$ with smooth boundary $\partial D$. The coefficient $\sigma(\cdot): \mathbb{R}\rightarrow \mathbb{R}$ is a deterministic continuous function. $B$ is a one dimensional standard Brownian motion defined on some filtrated probability space $(\Omega, {\cal F}, {\cal F}_t, P)$.  Throughout this paper, we assume that the initial value $u_0$ is a deterministic function in $L^2(D)$, and the system (\ref{1.a}) is considered in the space $L^2(D)$.

\vskip 0.3cm

Such a logarithmic nonlinearity has been introduced in the study of nonlinear wave mechanics and relativistic field in Physics \cite{R,BM}. The logarithmic wave mechanics and logarithmic Schr\"{o}dinger equations have been studied by many authors. We refer the readers to \cite{R,BM} and references therein  for details.
\vskip 0.3cm
Deterministic heat equations with logarithmic nonlinearity have been studied by several people.
Chen, Luo and Liu  in \cite{CLL} considered the deterministic heat equation with logarithmic nonlinearity on a bounded domain (see also \cite{CT}). They obtained the existence of global solutions and studied the blow up problem at infinity when the initial value $u_0\in L^2(D)$ satisfies certain energy  conditions. However, they didn't show the uniqueness of the solutions. In \cite{AC}, Alfaro and Carles considered the deterministic heat equation with logarithmic nonlinearity on the whole real line $\mathbb{R}$. They obtained the existence and uniqueness of classical solutions for a class of initial data which are bounded and sufficiently smooth.
\vskip 0.3cm

For stochastic heat equations with logarithmic nonlinearity, we mention the paper \cite{DKZ} by Dalang, Khoshnevisan and Zhang. The stochastic heat equations driven by space time white noise were considered in \cite{DKZ}. One of their results  claims  that any $L^2$-valued  solution will not blow up. However, the  existence and uniqueness of $L^2$-valued solutions were not proved. There exists a vast amount of references on  stochastic partial differential equations (SPDEs). SPDEs with the  so called monotone or locally monotone coefficients/nonlinearity were considered by many people in the literature. However, coefficients with logarithmic nonlinearity do not fall into this category.   We refer the readers to the monographs \cite{DZ} and \cite{PR} for SPDEs in general.

\vskip 0.3cm
In this paper we establish the existence and uniqueness of solutions of the stochastic heat equation (\ref{1.a}) in the space $L^2(D)$ for all initial values $u_0\in L^2(D)$. Based on a new estimate of the difference of two logarithmic terms and a nonlinear type of Gronwall's inequalities, we prove the uniqueness of the solutions in the $L^2(D)$ space when the diffusion coefficient $\sigma$ satisfies a local Lipschitz condition. To obtain the existence of solutions, we use the Galerkin methods. We first establish  the existence of a probabilistic weak solution by showing  the tightness of the approximate solutions and identifying any their limit as the solution of the stochastic heat equation. The existence of probabilistic strong solutions  then follows by appealing to the Yamada Watanabe theorem. We have two results on the existence of the solutions. The first one is obtained under the sublinear growth condition on the diffusion coefficient $\sigma$. In this case, we also obtain a global moment estimate of the solution. The second is established under the superlinear growth condition on $\sigma$. However, we do not have the moment estimate for the latter.

\vskip 0.3cm
The rest of the paper is organized as follows. In Section 2, we recall the logarithmic Sobolev inequality, and present the  framework for (\ref{1.a}) and give our hypotheses. Section 3 is devoted to the proof of uniqueness of solutions. In Section 4, we  establish the well-posedness of the approximating finite-dimensional stochastic differential equations. In section 5, we establish the tightness of the approximate solutions and  prove the existence of solutions under the sublinear growth condition on the coefficient $\sigma$. The solution is obtained by taking the limit of the Galerkin approximations. We also provide a moment estimate for the solution. Section 6 is to obtain the existence of the solution under the superlinear growth condition of $\sigma$.  Section 7 is the Appendix, which contains two nonlinear types of Gronwall's inequalities used in Section 3 and Section 5.

\section{Preliminaries and hypotheses}\label{S:2}
\setcounter{equation}{0}

In this section we will set up the framework.
We first introduce the following standard spaces. Let $H:=L^2(D)$. The norm and the inner product of $H$ are  denoted by $\Vert \cdot\Vert$ and $(\cdot, \cdot)$, respectively. Denote by $V$ the Sobolev space $H_0^1(D)$, that is,  the completion of the set of smooth functions with compact support ($C^{\infty}_0(D)$) under the norm
\begin{align}
  \Vert u\Vert_{V}^2=\int_D|\nabla u|^2(x)dx.
\end{align}
It is known that there exists an orthonormal basis $\{e_i\}_{i=1}^{\infty}$ of $H$ which consists of the eigenvectors of the negative Laplace operator under zero boundary conditions with corresponding eigenvalues $0<\lambda_i \uparrow<\infty$, that is
\begin{align}\label{A eigenvalue}
\Delta e_i=-\lambda_i e_i, \quad e_i|_{\partial D}=0,\quad i\in\mathbb{N}.
\end{align}
Moreover, $\{e_i\}_{i=1}^{\infty}$ is an orthogonal basis of $V$, and $e_i\in L^{\infty}(D)$.
Recall the $\rm Poincar\acute{e}$ inequality, i.e.
\begin{align}\label{Poincare}
\Vert u\Vert^2 \leq \frac{1}{\lambda_1}\Vert u\Vert_{V}^2, \quad \forall\, u\in V .
\end{align}
To handle the logarithmic term, we need the logarithmic Sobolev inequality of Gross \cite{G} in the following form.
For any $\varepsilon>0$ and $u\in V$, we have
\begin{align}\label{log sobolev}
\int_D |u(x)|^2\log|u(x)| dx \leq \varepsilon \Vert u \Vert_{V}^2 + \left(\frac{d}{4}\log\frac{1}{\varepsilon}\right) \Vert u\Vert^2 + \Vert u\Vert^2\log\Vert u\Vert.
\end{align}
Set
\[
\log_{+}z:=\log (1\vee z) .
\]
From the above logarithmic Sobolev inequality, it follows that for any $\varepsilon>0$ and $u\in V$,
\begin{align}\label{log sobolev modi}
& \int_D |u(x)|^2\log_{+}|u(x)| dx \nonumber\\
= & \int_D |u(x)|^2\log |u(x)| dx + \int_D |u(x)|^2\log \frac{1}{|u(x)|} \mathbf{1}_{\{0\leq |u|\leq 1\}}  dx \nonumber\\
\leq & \int_D |u(x)|^2\log |u(x)| dx + \frac{1}{2\mathrm{e}}m(D) \nonumber\\
\leq & \varepsilon \Vert u \Vert_{V}^2 + \left(\frac{d}{4}\log\frac{1}{\varepsilon}\right) \Vert u\Vert^2 + \Vert u\Vert^2\log\Vert u\Vert + \frac{1}{2\mathrm{e}}m(D) ,
\end{align}
where we have used
\begin{align}
  \max_{0\leq z\leq 1} z^2\log \frac{1}{z} =\frac{1}{2\mathrm{e}} .
\end{align}

If we identify the Hilbert space $H$ with its dual space $H^*$ by the Riesz representation, then we obtain a Gelfand triple
\begin{align*}
V\subset H\subset V^*.
\end{align*}

\noindent We denote by $\langle f,v\rangle$ the dual pairing between $f\in V^*$ and $v\in V$. It is easy to see that
\begin{align}\label{eq P4 star}
(u,v)=\langle u,v\rangle,\ \ \ \forall\,u\in H,\ \ \forall\,v\in V.
\end{align}

Set
\[\quad  u(t)(x):=u(t,x), \quad\quad \left(u(t)\log|u(t)|\right)(x):= u(t,x)\log|u(t,x)|. \]
Then (\ref{1.a}) can be formulated as the following stochastic evolution  equation
\begin{numcases}{}
  u(t) = u_0 +\int_0^t \Delta u(s)ds + \int_0^t u(s)\log|u(s)| ds + \int_0^t \sigma(u(s))dB_s , \nonumber\\
\label{Abstract}  u(0)=u_0\in H.
\end{numcases}

\begin{definition}
An $H$-valued continuous and $\mathcal{F}_t$-adapted stochastic process $u$ is called a  solution of (\ref{Abstract}), if the following two conditions hold:
\begin{itemize}
  \item [(i)] $u\in L^2([0,T];V)$ for any $T>0$, $P$\text{-a.s.}
  \item [(ii)] $u$ satisfies the equation (\ref{Abstract}) in $V^*$, $P$-a.s. for any $t\geq 0$.
\end{itemize}
\end{definition}

\vskip 0.5cm

Now we introduce our hypotheses on the diffusion coefficient $\sigma$. For the uniqueness, we assume that $\sigma$ satisfies a local Lipschitz condition.
\begin{itemize}
  \item [\hypertarget{H.1}{{\bf (H.1)}}] There exist constants $L_1$ and $L_2$  such that for all $x,y\in \mathbb{R}$,
\begin{gather}
\label{Local-Lipschitz} |\sigma(x)-\sigma(y)| \leq  L_1|x-y|+ L_2 |x-y| \left(\log_{+}(|x|\vee |y|)\right)^{\frac{1}{2}} .
\end{gather}
\end{itemize}
For the existence, we introduce two different hypotheses,  which lead to two different results.
\begin{itemize}
  \item [\hypertarget{H.2}{{\bf (H.2)}}] There exist constants $\theta\in[0,1)$, $C_1$ and $C_2$ such that for all $x\in \mathbb{R}$,
\begin{gather}
\label{sublinear} |\sigma(x)| \leq  C_1 + C_2 |x|^{\theta}.
\end{gather}
\end{itemize}

\begin{itemize}
  \item [\hypertarget{H.3}{{\bf (H.3)}}] There exist constants $C_3$ and $C_4$  such that for all $x\in \mathbb{R}$,
\begin{gather}
\label{superlinear} |\sigma(x)| \leq  C_3 + C_4 |x|\left(\log_{+} |x|\right)^{\frac{1}{2}}.
\end{gather}
\end{itemize}

\begin{remark}
  In particular, (\hyperlink{H.2}{H.2}) implies the linear growth of $\sigma$ and (\hyperlink{H.3}{H.3}). Obviously, (\hyperlink{H.1}{H.1}) implies (\hyperlink{H.3}{H.3}).
A typical example of functions satisfying (\hyperlink{H.1}{H.1}) and (\hyperlink{H.3}{H.3}) is
\begin{align}
\sigma(x)=
\begin{cases}
x\left(\log|x|\right)^{\frac{1}{2}} & |x|\geq \mathrm{e},\\
x & |x|\leq \mathrm{e} .
\end{cases}
\end{align}
\end{remark}


\section{Uniqueness of solutions}
\setcounter{equation}{0}

In this section, we will show the pathwise uniqueness of solutions to equation (\ref{Abstract}).  To do this, we first provide an estimate concerning the difference of two logarithmic terms.


\begin{lemma}\label{lemma 3.1}
For any $u,v\in V$, $\varepsilon>0$, and $\alpha\in(0,1)$, we have
\begin{align}\label{3.1}
  & \left(u\log|u| - v\log|v|, u-v\right) \nonumber\\
  \leq &  \varepsilon \Vert u-v\Vert_{V}^2 + \left(1+\frac{d}{4}\log \frac{1}{\varepsilon} \right)  \Vert u-v\Vert^2
   + \Vert u-v \Vert^2\log\Vert u-v\Vert \nonumber\\  & + \frac{1}{2(1-\alpha)\mathrm{e}}\left(\Vert u\Vert^{2(1-\alpha)}+\Vert v\Vert^{2(1-\alpha)}\right)\Vert u-v\Vert^{2\alpha}.
\end{align}
\end{lemma}

\vskip 0.3cm
\noindent {\bf Proof}.
Splitting the domain of integration we have
\begin{align}\label{3.2}
  & \left(u\log|u| - v\log|v|, u-v\right) \nonumber\\
   = & \int_{\{|u|>|v|\}} (u-v)^2\log|u| dx  +\int_{\{|u|>|v|\}}  v(\log|u|-\log|v|)(u-v)dx \nonumber\\
   +& \int_{\{|u|<|v|\}} (u-v)^2\log|v| dx  +\int_{\{|u|<|v|\}}  u(\log|u|-\log|v|)(u-v)dx \nonumber\\
   \leq & \int_D \big(u-v\big)^2\log(|u|\vee|v|) dx + \int_D |u-v|^2 dx,
\end{align}
where we have used the fact that for $0<x_1<x_2$,
\[ |\log(x_1)-\log(x_2)|\leq \frac{1}{x_1}|x_1-x_2|. \]
We now estimate the first term on the right of (\ref{3.2}),
\begin{align}\label{3.3}
   & \int_D \big(u-v\big)^2\log\big(|u|\vee|v|\big) dx \nonumber\\
 \leq & \int_D \big(u-v\big)^2\log|u- v| dx
   + \int_D \big(u-v\big)^2\log\left(\frac{|u|\vee|v|}{|u- v|}\right) dx .
\end{align}
The first term on the right of (\ref{3.3}) can be bounded by using the logarithmic  Sobolev inequality (see (\ref{log sobolev})):
\begin{align}\label{3.4}
  & \int_D \big(u-v\big)^2\log|u- v| dx \nonumber \\
  \leq & \varepsilon\Vert u-v\Vert_{V}^2 + \left(\frac{d}{4}\log\frac{1}{\varepsilon}\right)  \Vert u-v\Vert^2   + \Vert u-v \Vert^2\log\Vert u-v\Vert  .
\end{align}
For the second term on the right of (\ref{3.3}), we have
\begin{align}\label{3.5}
   \int_D (u-v)^2\log\left(\frac{|u|\vee|v|}{|u- v|}\right) dx
  = & \int_{\{uv\leq 0\}}(u-v)^2\log\left(\frac{|u|\vee|v|}{|u- v|}\right) dx \nonumber\\
  & + \int_{\{uv> 0, |u|\geq |v|\}} (u-v)^2\log\left(\frac{u}{u-v}\right) dx \nonumber\\
  & + \int_{\{uv> 0, |u|<|v|\}}(u-v)^2\log\left(\frac{v}{v-u}\right) dx \nonumber\\
  =: &\  I_1 +I_2 +I_3 .
\end{align}
Note that
\begin{align}\label{3.6}
  I_1\leq 0 .
\end{align}
For any fixed $b>0$, $p>0$, we define a function $h:[0,b]\rightarrow\mathbb{R}$ by
\[h(z):=(b-z)^{p}\log\left(\frac{b}{b-z}\right) .\]
A straightforward calculation leads to
\[\max_{0\leq z\leq b}h(z)= \frac{1}{p\mathrm{e}}b^p .\]
Therefore, by H\"o{}lder's inequality we get
\begin{align}\label{3.7}
  I_2\leq & \int_{\{uv> 0, |u|\geq |v|\}} |u-v|^{2\alpha}|u-v|^{2(1-\alpha)}\log\left(\frac{u}{u-v}\right) dx \nonumber \\
  \leq & \int_{\{uv> 0, |u|\geq |v|\}} |u-v|^{2\alpha}\times\frac{1}{2(1-\alpha)\mathrm{e}}|u|^{2(1-\alpha)} dx \nonumber \\
  \leq & \frac{1}{2(1-\alpha)\mathrm{e}} \left[\int_D |u-v|^{2}dx\right]^{\alpha} \times  \left[\int_D |u|^{2}dx\right]^{1-\alpha} \nonumber \\
  = & \frac{1}{2(1-\alpha)\mathrm{e}} \Vert u\Vert^{2(1-\alpha)} \Vert u-v\Vert^{2\alpha} .
\end{align}
Similarly,
\begin{align}\label{3.8}
  I_3\leq \frac{1}{2(1-\alpha)\mathrm{e}} \Vert v\Vert^{2(1-\alpha)} \Vert u-v\Vert^{2\alpha}.
\end{align}
Putting (\ref{3.2})-(\ref{3.8}) together, we obtain (\ref{3.1}).
$\blacksquare$

\begin{lemma}\label{lemma20190628.1}
  For any $u,v\in V$, $\epsilon>0$, and $\alpha\in(0,1)$, we have
\begin{align}\label{20190629.1433.1}
  & \int_D |u(x)-v(x)|^2 \log_{+}\left(|u(x)|\vee |v(x)|\right) dx \nonumber\\
  \leq &  \epsilon \Vert u-v\Vert_{V}^2 + \left(\frac{d}{4}\log \frac{1}{\epsilon} \right)  \Vert u-v\Vert^2
   + \Vert u-v \Vert^2\log\Vert u-v\Vert \nonumber\\  & + \frac{1}{2(1-\alpha)\mathrm{e}}\left(\Vert u\Vert^{2(1-\alpha)}+\Vert v\Vert^{2(1-\alpha)}\right)\Vert u-v\Vert^{2\alpha} \nonumber\\
   & + \frac{1}{2(1-\alpha)\mathrm{e}} \left(4m(D)\right)^{1-\alpha} \Vert u-v\Vert^{2\alpha},
\end{align}
where $m(D)$ is the Lebesgue measure of domain $D$.
\end{lemma}
\noindent {\bf Proof}. Note that
\begin{align}\label{190629.1916.1}
  & \int_D |u(x)-v(x)|^2 \log_{+}\left(|u(x)|\vee |v(x)|\right) dx \nonumber\\
  = & \int_D |u(x)-v(x)|^2 \log \left(|u(x)|\vee |v(x)|\right) dx \nonumber\\
  & + \int_D |u(x)-v(x)|^2 \log \frac{1}{|u(x)|\vee |v(x)|}  \mathbf{1}_{\{|u(x)|\vee |v(x)|\leq 1\}} dx  \nonumber\\
  =& :  J_1 +J_2 .
\end{align}
The terms $J_1$ and $J_2$ are estimated as follows. On the one hand, combining (\ref{3.3})-(\ref{3.8}) together, we have
\begin{align}\label{190629.1916.2}
  J_1 \leq & \epsilon \Vert u-v\Vert_{V}^2 + \left(\frac{d}{4}\log \frac{1}{\epsilon} \right)  \Vert u-v\Vert^2
   + \Vert u-v \Vert^2\log\Vert u-v\Vert \nonumber\\  & + \frac{1}{2(1-\alpha)\mathrm{e}}\left(\Vert u\Vert^{2(1-\alpha)}+\Vert v\Vert^{2(1-\alpha)}\right)\Vert u-v\Vert^{2\alpha}  .
\end{align}
On the other hand,
\begin{align}\label{190629.1916.3}
 J_2
  = & \int_{|v|\leq |u|} |u(x)-v(x)|^{2\alpha} |u(x)-v(x)|^{2(1-\alpha)}\log \frac{1}{|u(x)|}  \mathbf{1}_{\{|v(x)|\leq |u(x)|\leq 1\}} dx \nonumber\\
  & + \int_{|u|\leq |v|} |u(x)-v(x)|^{2\alpha} |u(x)-v(x)|^{2(1-\alpha)}\log \frac{1}{|v(x)|}  \mathbf{1}_{\{|u(x)|\leq |v(x)|\leq 1\}} dx \nonumber\\
\leq & \int_{|v|\leq |u|} |u(x)-v(x)|^{2\alpha} (2|u(x)|)^{2(1-\alpha)}\log \frac{1}{|u(x)|}  \mathbf{1}_{\{|u(x)|\leq 1\}} dx \nonumber\\
  & + \int_{|u|\leq |v|} |u(x)-v(x)|^{2\alpha} (2|v(x)|)^{2(1-\alpha)}\log \frac{1}{|v(x)|}  \mathbf{1}_{\{|v(x)|\leq 1\}} dx \nonumber\\
  \leq & \frac{2^{2(1-\alpha)}}{2(1-\alpha)\mathrm{e}} \int_{D} |u(x)-v(x)|^{2\alpha} \mathbf{1}_{\{|u(x)|\vee |v(x)|\leq 1\}} dx \nonumber\\
  \leq & \frac{4^{1-\alpha}}{2(1-\alpha)\mathrm{e}} \Vert u-v\Vert^{2\alpha} (m(D))^{1-\alpha} ,
\end{align}
where we have used
\begin{align}
  \max_{0\leq z\leq 1} z^{2(1-\alpha)} \log\frac{1}{z} = \frac{1}{2(1-\alpha)\mathrm{e}} .
\end{align}
Combining (\ref{190629.1916.1})-(\ref{190629.1916.3}) together, we obtain (\ref{20190629.1433.1}).
$\blacksquare$

\begin{theorem}\label{thm 3.2}
Suppose hypothesis (\hyperlink{H.1}{H.1}) holds. Then the pathwise uniqueness holds for equation (\ref{Abstract}) in $L^2(D)$.
\end{theorem}
\noindent {\bf Proof}.
Let $u,v$ be two solutions of equation (\ref{Abstract}). For any $M>0$ and $1\geq \delta>0$, we define stopping times
\begin{align*}
  \tau_M:=& \inf\{t>0: \Vert u(t)\Vert^2\vee\Vert v(t)\Vert^2>M\}, \\
  \tau_M^{\prime}:=& \inf\{t>0: \int_0^t\Vert u(s)\Vert_{V}^2ds>M\}, \\
   \tau^{\delta}:=& \inf\{t>0: \Vert u(t)-v(t)\Vert>\delta\} ,\\
\tau_M^{\delta}:=& \tau_M\wedge \tau_M^{\prime} \wedge \tau^{\delta}.
\end{align*}
Set
\[Z(t):=u(t)-v(t) .\]
Applying Ito's formula, we have
\begin{align}\label{3.9}
  & \Vert Z(t\wedge\tau_M^{\delta})\Vert^2 + 2\int_0^{t\wedge\tau_M^{\delta}}\Vert Z(s)\Vert_{V}^2 ds \nonumber\\
  = & 2\int_0^{t\wedge\tau_M^{\delta}}\big(u(s)\log|u(s)|-v(s)\log|v(s)|, u(s)-v(s)\big)ds \nonumber\\
  & +  2\int_0^{t\wedge\tau_M^{\delta}}\big(\sigma(u(s))-\sigma(v(s)),u(s)-v(s)\big)dB_s \nonumber\\
  & +  \int_0^{t\wedge\tau_M^{\delta}}\Vert\sigma(u(s))-\sigma(v(s))\Vert^2 ds \nonumber\\
  \leq & 2\int_0^{t\wedge\tau_M^{\delta}}\big(u(s)\log|u(s)|-v(s)\log|v(s)|, u(s)-v(s)\big)ds \nonumber\\
  & +  2\int_0^{t\wedge\tau_M^{\delta}}\big(\sigma(u(s))-\sigma(v(s)),u(s)-v(s)\big)dB_s \nonumber\\
  & + 2L_1^2 \int_0^{t\wedge\tau_M^{\delta}}\Vert Z(s)\Vert^2 ds \nonumber\\
  & + 2L_2^2\int_0^{t\wedge\tau_M^{\delta}}\int_D |u(s,x)-v(s,x)|^2 \log_{+}\left(|u(s,x)|\vee |v(s,x)|\right) dx ds ,
\end{align}
where we have used hypothesis (\hyperlink{H.1}{H.1}) in the last inequality.
Substituting estimate (\ref{3.1}) with $\varepsilon=\frac{1}{4}$ and estimate (\ref{20190629.1433.1}) with $\epsilon=\frac{1}{4L_2^2}$  into the above equality, we get
\begin{align}\label{3.10}
  & \Vert Z(t\wedge\tau_M^{\delta})\Vert^2 + \int_0^{t\wedge\tau_M^{\delta}}\Vert Z(s)\Vert_{V}^2 ds \nonumber\\
  \leq & C\int_0^{t\wedge\tau_M^{\delta}}\Vert Z(s)\Vert^2 ds \nonumber\\
  & + (2+2L_2^2)\int_0^{t\wedge\tau_M^{\delta}}\Vert Z(s)\Vert^2 \log \Vert Z(s)\Vert ds \nonumber\\
  & + \frac{1+L_2^2}{(1-\alpha)\mathrm{e}}\int_0^{t\wedge\tau_M^{\delta}} \left(\Vert u(s)\Vert^{2(1-\alpha)} + \Vert v(s)\Vert^{2(1-\alpha)}\right) \Vert Z(s)\Vert^{2\alpha} \nonumber\\
  & + \frac{L_2^2}{(1-\alpha)\mathrm{e}}\int_0^{t\wedge\tau_M^{\delta}} \left(4m(D)\right)^{1-\alpha} \Vert Z(s)\Vert^{2\alpha} \nonumber\\
  & +  2\int_0^{t\wedge\tau_M^{\delta}}\big(\sigma(u(s))-\sigma(v(s)),u(s)-v(s)\big)dB_s .
\end{align}
By the definition of $\tau_M^{\delta}$,
\begin{align}
\int_0^{t\wedge\tau_M^{\delta}}\Vert Z(s)\Vert^2 \log \Vert Z(s)\Vert ds \leq 0 .
\end{align}
Hence taking expectations on both sides of (\ref{3.10})
and setting
\[Y(t):=E\Vert Z(t\wedge\tau_M^{\delta})\Vert^2 \]
yields
\begin{align}\label{3.11}
  Y(t)\leq C\int_0^t Y(s) ds + \frac{2(1+L_2^2)M^{1-\alpha}+L_2^2(4m(D))^{1-\alpha}}{(1-\alpha)\mathrm{e}}\int_0^t Y(s)^{\alpha} ds .
\end{align}
From Lemma \ref{A.1} in Appendix, it follows that
\begin{align}\label{3.12}
  Y(t)\leq & \left\{(1-\alpha)\int_{0}^t \frac{2(1+L_2^2)M^{1-\alpha}+L_2^2(4m(D))^{1-\alpha}}{(1-\alpha)\mathrm{e}}\times e^{(1-\alpha)\times C(t-s)}ds\right\}^{\frac{1}{1-\alpha}} \nonumber\\
  = & \left[\frac{2(1+L_2^2)M^{1-\alpha}+L_2^2(4m(D))^{1-\alpha}}{\mathrm{e}}\right]^{\frac{1}{1-\alpha}} \times \left(\int_{0}^{t} e^{(1-\alpha)\times Cs }ds\right)^{\frac{1}{1-\alpha}} \nonumber\\
  \leq & \frac{1}{2} \left\{\left[\frac{4(1+L_2^2)t^{\alpha}}{\mathrm{e}}\right]^{\frac{1}{1-\alpha}}\times M + \left[ \frac{2L_2^2 t^{\alpha}}{\mathrm{e}}\right]^{\frac{1}{1-\alpha}} \times 4m(D)\right\} \times \left(\int_{0}^{t}e^{Cs}ds\right) . \nonumber\\
\end{align}
Letting $\alpha\rightarrow 1$ and setting $T^*:=\left(\frac{\mathrm{e}}{4(1+L_2^2)}\right)^2$, we obtain
\begin{align}
Y(t)=0, \quad \forall\, 0\leq t\leq T^* .
\end{align}
Since this time interval is independent of  the initial value, starting from time $T^*$, by the same argument, it can be deduced that $Y(t)=0$ for any $t\in [T^*, 2T^*]$. Repeating this argument, we deduce that $Y(t)=0$ for any $t\geq 0$. This means
\begin{align}\label{3.13}
  E\Vert Z(t\wedge\tau_M\wedge\tau_M^{\prime}\wedge\tau^{\delta})\Vert^2=0 , \quad \forall\,t\geq 0.
\end{align}
Since $u, v$ are two global solutions, $\tau_M\rightarrow\infty$ and $\tau_M^{\prime}\rightarrow\infty$,  $P$-a.s. as $M\rightarrow\infty$. Letting $M\rightarrow\infty$ in (\ref{3.13}), we get
\begin{align}
  E\Vert Z(t\wedge\tau^{\delta})\Vert^2=0 , \quad \forall\,t\geq 0.
\end{align}
This implies
\begin{align}
  P(\tau^{\delta}>t) =1 , \quad \forall\,t\geq 0, \forall\, \delta>0.
\end{align}
Hence
\begin{align}
u(t)=v(t), \quad P\text{-a.s.}, \quad \forall\,t\geq 0.
\end{align}
The pathwise uniqueness follows from the path continuity of $u, v$ in $H$.
$\blacksquare$

\section{Galerkin approximating solutions}
\setcounter{equation}{0}

We will employ the Galerkin methods to prove the existence of solutions. To this end,
we first study the well-posedness of the Galerkin approximating equations in this section.
\vskip 0.3cm
Let $H_n$ denote the $n$-dimensional subspace of $H$ spanned by $\{e_1,\dots, e_n\}$. Let $P_n: V^*\rightarrow H_n$ be defined by
\begin{align}\label{4.1}
  P_n g := \sum_{i=1}^{n}\langle g,e_i\rangle e_i .
\end{align}
For any integer $n\geq 1$, we consider the following stochastic differential equation in the finite-dimensional space  $H_n$:
\begin{numcases}{}
  du_n(t)= \Delta u_n(t)dt+P_n [u_n(t)\log|u^n(t)|]dt+ P_n\sigma(u_n(t))dB_t, \quad t\geq 0, \nonumber\\
\label{4.2}  u_n(0)= P_n u_0 ,
\end{numcases}
such that
\begin{align}\label{4.3}
  u_n(t)=\sum_{i=1}^{n}g_{in}(t)e_i .
\end{align}
$u_n$ solves (\ref{4.2}) if and only if $\{g_{in}\}_{i=1}^{n}$ solves the system
\begin{align}\label{4.4}
   dg_{jn}(t)= & d(u^n(t),e_i) \nonumber\\
  =& -\lambda_j g_{jn}(t)dt + \left(\sum_{i=1}^n g_{in}(t)e_i\log|\sum_{i=1}^n g_{in}(t)e_i|, e_j\right) dt \nonumber\\
   & + \left(\sigma\Big(\sum_{i=1}^n g_{in}(t)e_i\Big), e_j\right)dB_t, \quad j=1,2 ,\dots, n.
\end{align}

To present results on the existence and uniqueness of equation
(\ref{4.4}), we introduce the following  functions $F_j$ and $G_j$, $j=1,\dots, n$,  on $\mathbb{R}^n$,
\begin{align}\label{4.5}
  F_{j}(y_1,\dots, y_n):&=\int_D e_j(x)\left(\sum_{i=1}^n y_i e_i(x)\right)\log\left|\sum_{i=1}^n y_i e_i(x)\right| dx  , \nonumber\\
  G_{j}(y_1,\dots, y_n):&=\int_D e_j(x)\sigma\left(\sum_{i=1}^n y_i e_i(x)\right) dx .
\end{align}
In the following, the length of a vector $y=(y_1,\dots, y_n)\in\mathbb{R}^n$ is denoted by $|y|$. Here are some estimates for the function $F_j$, $j=1,\dots, n$.
\begin{lemma}\label{lemma 4.1}

\begin{itemize}
  \item [(i)] There exist constants $\widetilde{L_1}$, $\widetilde{L_2}$, $\widetilde{L_3}$ and $\delta>0$ such that for any $y, z\in\mathbb{R}^n$ and $|y-z|\leq \delta$,
\begin{align}
    & \left|F_{j}(y_1,\dots, y_n)-F_{j}(z_1,\dots, z_n)\right| \nonumber\\
     \leq & \widetilde{L_1}|y-z| + \widetilde{L_2}|y-z|\log_{+}(|y|\vee |z|) + \widetilde{L_3}|y-z|\log\frac{1}{|y-z|} .
\end{align}
  \item [(ii)] There exist constants $\widetilde{C_1}$, $\widetilde{C_2}$  such that for any $y\in \mathbb{R}^n$,
  \begin{align}\label{190629.2235.1}
    \left|F_{j}(y_1,\dots, y_n)\right| \leq  \widetilde{C_1} + \widetilde{C_2}|y|\log_{+}|y| .
  \end{align}
\end{itemize}

\end{lemma}
\begin{remark}
   (i) of Lemma \ref{lemma 4.1} implies that $F_j$, $j=1,\dots, n$, satisfy the local $\log$-Lipschitz condition, i.e. there exists a constant $0<\delta<1$, and for any $r>0$, there exists a constant $C_r>0$ such that
   \begin{align}\label{4.7}
     \left|F_{j}(y_1,\dots, y_n)-F_{j}(z_1,\dots, z_n)\right|\leq C_r |y-z|\log\frac{1}{|y-z|} ,
   \end{align}
   for any $|y|\vee|z|\leq r$ and $|y-z|< \delta$.
\end{remark}
\noindent {\bf Proof}. Proof of (i). Take
\begin{align}
  \delta := \min\left\{\left(\sum_{i=1}^n \Vert e_i\Vert_{L^{\infty}}^2\right)^{-\frac{1}{2}} , \frac{\sqrt{m(D)}}{\mathrm{e}} \right\},
\end{align}
where $m(D)$ is the Lebesgue measure of domain $D$.
For simplicity, we introduce two functions $v_1$ and $v_2$ as
\[ v_1(x):=\sum_{i=1}^n y_i e_i(x),\quad v_2(x):=\sum_{i=1}^n z_i e_i(x) .\]
Then
\begin{align}\label{4.8}
  & \left|F_{j}(y_1,\dots, y_n)-F_{j}(z_1,\dots, z_n)\right| \nonumber\\
  = & \bigg|\int_{\{|v_1|>|v_2|\}}e_j (v_1 -v_2)\log|v_1| dx + \int_{\{|v_1|>|v_2|\}}e_j v_2\left(\log|v_1|-\log|v_2|\right) dx \nonumber\\
  & +\int_{\{|v_1|<|v_2|\}}e_j (v_1 -v_2)\log|v_2| dx + \int_{\{|v_1|<|v_2|\}}e_j v_1\left(\log|v_1|-\log|v_2|\right) dx \bigg| \nonumber\\
  \leq & \bigg|\int_D e_j(v_1 -v_2)\log(|v_1|\vee|v_2|)dx\bigg| + \int_D |e_j| |v_1 -v_2|dx \nonumber\\
  \leq & \Vert e_j\Vert_{L^{\infty}}\times\int_D |v_1 -v_2|\left|\log(|v_1|\vee|v_2|)\right|dx + |y-z|  .
\end{align}
%
%
%

\noindent Note that for any $x\in D$,
\begin{align}\label{4.9}
  |v_1 (x)|=\left|\sum_{i=1}^n y_i e_i(x)\right|\leq & |y|\times \left(\sum_{i=1}^n \Vert e_i\Vert_{L^{\infty}}^2\right)^\frac{1}{2} \nonumber\\
  \leq & (|y|\vee |z|)\times\left(\sum_{i=1}^n \Vert e_i\Vert_{L^{\infty}}^2\right)^\frac{1}{2} ,
\end{align}
and
\begin{align}\label{4.10}
  \log\frac{1}{|a|\vee|b|}\leq \log\frac{1}{|a-b|} +\log 2 , \quad \forall\, a,b\in [-1,1] .
\end{align}
By (\ref{4.9}) and (\ref{4.10}), we get
\begin{align}\label{4.11}
  &\int_D |v_1 -v_2|\left|\log(|v_1|\vee|v_2|)\right|dx \nonumber\\
  = & \int_{\{1< |v_1|\vee|v_2|\}} |v_1 -v_2||\log(|v_1|\vee |v_2|)| dx \nonumber\\
  & + \int_{\{1\geq |v_1|\vee|v_2|\}} |v_1 -v_2||\log(|v_1|\vee |v_2|)| dx \nonumber\\
  \leq & \left[\log_{+}\left((|y|\vee|z|)\times\Big(\sum_{i=1}^n \Vert e_i\Vert_{L^{\infty}}^2\Big)^{\frac{1}{2}}\right)\right]\times\int_D |v_1 -v_2| dx \nonumber\\
  & + \int_{\{1\geq |v_1|\vee|v_2|\}} |v_1 -v_2|\left(\log\frac{1}{|v_1-v_2|}+\log 2\right) dx \nonumber\\
  = & \left[\log_{+}\left((|y|\vee|z|)\times\Big(\sum_{i=1}^n \Vert e_i\Vert_{L^{\infty}}^2\Big)^{\frac{1}{2}}\right) +\log 2\right] \times\int_D |v_1 -v_2| dx \nonumber\\
  & +\int_{\{1\geq |v_1|\vee|v_2|\}} |v_1 -v_2| \log\frac{1}{|v_1-v_2|}  dx \nonumber\\
  =& : I_1 + I_2 .
\end{align}
Similar to (\ref{4.9}) and using $|y-z|\leq \delta$, we have for any $x\in D$,
\begin{align}
  |v_1 (x) -v_2 (x)| \leq |y-z| \left(\sum_{i=1}^n \Vert e_i\Vert_{L^{\infty}}^2\right)^\frac{1}{2} \leq 1.
\end{align}
Thus
\begin{align}
  I_2 \leq \int_{D} |v_1 -v_2| \log\frac{1}{|v_1-v_2|}  dx .
\end{align}
Note that the function $x\mapsto x\log\frac{1}{x}$ is concave on $\mathbb{R}_+$. By Jensen's inequality, we have
\begin{align}\label{4.12}
  I_2  \leq  \int_D |v_1 -v_2| dx \times\log\frac{m(D)}{\int_D |v_1 -v_2| dx} .
\end{align}
Applying H\"o{}lder's inequality yields
\begin{align}\label{4.13}
  \int_D |v_1-v_2|dx \leq \sqrt{m(D)}\times \int_D |v_1-v_2|^2 dx = \sqrt{m(D)}\times |y-z| .
\end{align}
Hence the condition $|y-z|\leq \delta$ gives
\begin{align}\label{4.13.1}
  \int_D |v_1-v_2|dx  \leq \frac{m(D)}{\mathrm{e}} .
\end{align}
By (\ref{4.12})-(\ref{4.13.1}) and the fact that the function $x\mapsto x\log\frac{m(D)}{x}$ is increasing in $x\in [0,\frac{m(D)}{\mathrm{e}}]$, it can be seen that
\begin{align}\label{4.14}
  I_2\leq & \sqrt{m(D)}|y-z|\log\left(\frac{\sqrt{m(D)}}{|y-z|}\right) \nonumber\\
  \leq & \frac{1}{2}\sqrt{m(D)}\log \big(m(D)\big)|y-z| + \sqrt{m(D)}|y-z|\log\frac{1}{|y-z|} .
\end{align}
Putting (\ref{4.8}), (\ref{4.11}), (\ref{4.13}) and (\ref{4.14}) together, we obtain
\begin{align}\label{4.15}
  & \left|F_{j}(y_1,\dots, y_n)-F_{j}(z_1,\dots, z_n)\right| \nonumber\\
  \leq & \bigg[1+\Vert e_j\Vert_{L^{\infty}}\sqrt{m(D)}\times\bigg(\log_{+}\Big((|y|\vee|z|)\times\Big(\sum_{i=1}^n \Vert e_i \Vert_{L^{\infty}}^2\Big)^{\frac{1}{2}}\Big) +\log2 \nonumber\\
  & + \frac{1}{2}\log\big( m(D)\big)\bigg)\bigg] \times|y-z| + \Vert e_j\Vert_{L^{\infty}} \sqrt{m(D)} |y-z|\log\frac{1}{|y-z|} .
\end{align}
Therefore, there exist three constants $\widetilde{L_1}$, $\widetilde{L_2}$ and $\widetilde{L_3}$ such that
\begin{align}
    & \left|F_{j}(y_1,\dots, y_n)-F_{j}(z_1,\dots, z_n)\right| \nonumber\\
     \leq & \widetilde{L_1}|y-z| + \widetilde{L_2}|y-z|\log_{+}(|y|\vee |z|) + \widetilde{L_3}|y-z|\log\frac{1}{|y-z|} ,
\end{align}
for any $|y-z|\leq \delta$.

\vskip 0.5cm
Proof of (ii). Take $z=0$ in (\ref{4.8}), i.e. $v_2=0$. In this case, the term $I_2$ in (\ref{4.11}) can be bounded as
\begin{align}\label{190629.2232.1}
  I_2 \leq \int_{\{1>|v_1|\}} |v_1| \log\frac{1}{|v_1|} dx \leq \frac{m(D)}{\mathrm{e}},
\end{align}
where we have used
\begin{align}
  \max_{0\leq r\leq 1} r\log\frac{1}{r} \leq \frac{1}{\mathrm{e}}.
\end{align}
From (\ref{4.8}), (\ref{4.11}), (\ref{4.13}) and (\ref{190629.2232.1}), it follows that
\begin{align}
  & \left|F_{j}(y_1,\dots, y_n)\right| \nonumber\\
  \leq & \bigg[1+\Vert e_j\Vert_{L^{\infty}}\sqrt{m(D)}\times\bigg(\log_{+}\Big(|y|\times\Big(\sum_{i=1}^n \Vert e_i \Vert_{L^{\infty}}^2\Big)^{\frac{1}{2}}\Big) +\log2 \bigg)\bigg] \times|y| \nonumber\\
  & + \Vert e_j\Vert_{L^{\infty}} \times \frac{m(D)}{\mathrm{e}}.
\end{align}
Therefore, (\ref{190629.2235.1}) is obtained. This completes the proof of Lemma \ref{lemma 4.1}.
$\blacksquare$

\vskip 0.3cm

Using the similar methods as the proof of Lemma \ref{lemma 4.1}, it can be seen that the following estimates for the function $G_j$, $j=1,\dots, n$ hold.
\begin{lemma}\label{lemma 190629.2317}
\begin{itemize}
  \item [(i)] Suppose hypothesis (\hyperlink{H.1}{H.1}) holds. Then there exist constants $\widetilde{L_4}$ and $\widetilde{L_5}$ such that for any $y,z \in\mathbb{R}^n$,
\begin{align}
    & \left|G_{j}(y_1,\dots, y_n)-G_{j}(z_1,\dots, z_n)\right| \nonumber\\
     \leq & \widetilde{L_4}|y-z| + \widetilde{L_5}|y-z| \left(\log_{+}(|y|\vee |z|)\right)^{\frac{1}{2}}.
\end{align}
  \item [(ii)] Suppose hypothesis (\hyperlink{H.3}{H.3}) holds. Then there exist constants $\widetilde{C_3}$, $\widetilde{C_4}$  such that for any $y\in \mathbb{R}^n$,
  \begin{align}\label{190703.18585.1}
    \left|G_{j}(y_1,\dots, y_n)\right| \leq  \widetilde{C_3} + \widetilde{C_4}|y|\left(\log_{+}|y|\right)^{\frac{1}{2}} .
  \end{align}
\end{itemize}

\end{lemma}

\vskip 0.5cm

By Lemma \ref{lemma 4.1} and Lemma \ref{lemma 190629.2317}, we can directly apply  Theorem A, Theorem B and Theorem D in \cite{FZ} to obtain the following result.
\begin{theorem}\label{190630.0008}
If hypothesis (\hyperlink{H.3}{H.3}) holds, then stochastic differential equation (\ref{4.2}) has no explosion. If hypothesis (\hyperlink{H.1}{H.1}) holds, then there exists a unique global probabilistic strong solution to equation (\ref{4.2}).
\end{theorem}

\section{Existence of solutions: part I}
\setcounter{equation}{0}

In this part, we assume (\hyperlink{H.2}{H.2}) holds. From Theorem \ref{190630.0008}, the Galerkin approximating equation (\ref{4.2}) has a global solution. To prove the existence of solutions to (\ref{Abstract}), we will show the tightness of Galerkin approximating solutions. Passing to the limit, we first obtain the existence of probabilistic weak solutions in a short time interval; we then construct a global solution by piecing together the solutions over subintervals. We finally establish global moment estimates of the solution.
%

\subsection{Tightness of approximating solutions}

We are going to prove the tightness of the solutions $\{u_n, n\geq1\}$. To this end,
We first prepare some estimates.
\begin{lemma}\label{lemma 190706.1431}
Let
\begin{align}\label{Tp}
  T_{p}:=\log\frac{p}{p-1+\theta} ,
\end{align}
where $\theta$ is the constant appeared in hypothesis (\hyperlink{H.2}{H.2}).
Then $T_p$ is decreasing in $p\in[2,\infty)$. And moreover, under hypothesis (\hyperlink{H.2}{H.2}), we have for any $p\geq 2$,
\begin{align}\label{4.16}
  &\sup_{n}E\left[\sup_{s\in[0,T_p]}\Vert u_n(s)\Vert^p + \int_0^{T_p} \Vert u_n(s)\Vert^{p-2} \Vert u_n(s)\Vert_{V}^2ds  \right] \nonumber\\
  &\leq C_{p,\theta}\left(1+ \Vert u_0\Vert^{\frac{p^2}{p-1+\theta}}\right) <\infty,
\end{align}
for some constant $C_{p,\theta}$.
\end{lemma}
\noindent {\bf Proof}.
For any $n\in \mathbb{N}$, $M>0$, we define stopping times
\[
  \tau_M^n:=\inf\{t\geq 0: \Vert u_n(t)\Vert >M \}\wedge T_p.
\]
Since $u_n$ has no explosion, $\tau_M^n\uparrow T_p$, $P$-a.s. as $M\rightarrow\infty$. Applying Ito's formula, we have for $t\leq \tau_M^n$,
\begin{align}\label{4.17}
  d\Vert u_n(s)\Vert^2 = & - 2 \Vert \nabla u_n(s)\Vert^2ds + 2\big(u_n(s)\log|u_n(s)|,u_n(s)\big) ds  \nonumber\\
  & + \Vert P_n\sigma(u_n(s))\Vert^2 ds + 2 \big(\sigma(u_n(s)),u_n(s)\big) dB_s .
\end{align}
Once again applying Ito's formula gives
\begin{align}\label{4.18}
  \Vert u_n(t)\Vert^p = & \Vert u_0 \Vert^p - p\int_0^t \Vert u_n(s)\Vert^{p-2}\Vert u_n(s)\Vert_{V}^2 ds \nonumber\\
  & + p\int_0^t \Vert u_n(s)\Vert^{p-2} \big(u_n(s)\log|u_n(s)|,u_n(s)\big) ds \nonumber\\
  & + \frac{p}{2} \int_0^t \Vert u_n(s)\Vert^{p-2}\Vert P_n\sigma(u(s))\Vert^2 ds \nonumber\\
  & + p\int_0^t  \Vert u_n(s)\Vert^{p-2}\big(\sigma(u_n(s)),u_n(s)\big)dB_s \nonumber\\
  & + \frac{p(p-2)}{2} \int_0^t \Vert u_n(s)\Vert^{p-4}\big(\sigma(u_n(s)),u_n(s)\big)^2 ds .
\end{align}
An appeal to the logarithmic Sobolev inequality (\ref{log sobolev}) with $\varepsilon=1/2$ yields
\begin{align}\label{4.19}
    & \Vert u_n(t)\Vert^p + p\int_0^t \Vert u_n(s)\Vert^{p-2}\Vert u_n(s)\Vert_{V}^2 ds \nonumber\\
  \leq  & \Vert u_0 \Vert^p + p\int_0^t \Vert u_n(s)\Vert^{p-2} \bigg(\frac{1}{2}\Vert u_n(s)\Vert_{V}^2+\frac{d\log 2}{4}\Vert u_n(s)\Vert^2 \nonumber \\
  & + \Vert u(s)\Vert^2\log\Vert u_n(s)\Vert\bigg) ds  + p\int_0^t  \Vert u_n(s)\Vert^{p-2}\big(\sigma(u_n(s)),u_n(s)\big)dB_s \nonumber\\
  & + \frac{p(p-1)}{2} \int_0^t \Vert u_n(s)\Vert^{p-2}\Vert\sigma(u_n(s))\Vert^2 ds .
\end{align}
Under hypothesis (\hyperlink{H.2}{H.2}), $\sigma$ also satisfies the linear growth condition, so we have
\begin{align}\label{4.20}
      & \Vert u_n(t)\Vert^p + \frac{p}{2}\int_0^t \Vert u_n(s)\Vert^{p-2}\Vert u_n(s)\Vert_{V}^2 ds \nonumber\\
  \leq  & M(t)+ C\int_0^t \Vert u_n(s)\Vert^p ds + \int_0^t \Vert u_n(s)\Vert^{p}\log\Vert u_n(s)\Vert^p ds ,
\end{align}
where
\begin{align}\label{4.21}
  M(t):=\Vert u_0 \Vert^p + C + p\sup_{r\in [0,t]}\left|\int_0^r  \Vert u_n(s)\Vert^{p-2}\big(\sigma(u_n(s)),u_n(s)\big)dB_s\right| ,
\end{align}
and these  constants $C$ are independent of $n$.
Applying the log-Gronwall inequality (see Lemma \ref{A.2} in Appendix) to (\ref{4.20}), we get for $t\in[0,\tau_M^n]$,
\begin{align}\label{4.22}
  \Vert u_n(t)\Vert^p + \frac{p}{2}\int_0^t \Vert u_n(s)\Vert^{p-2}\Vert u_n(s)\Vert_{V}^2 ds \leq \left(1\vee M(t) \right)^{e^{t}}\times e^{C (e^t -1)}.
\end{align}
Set
\[
  X_n(t):=E\left[\sup_{s\in[0,t\wedge\tau_M^n]}\Vert u_n(s)\Vert^p\right] .
\]
Note that
\begin{align}\label{4.23}
  \Vert \sigma(u)\Vert    \leq C \left(\int_D (1+|u(x)|^{\theta})^2 dx\right)^\frac{1}{2}
   \leq C(1+\Vert u\Vert^{\theta}) .
\end{align}
By BDG inequality and Young's inequality $ab\leq \epsilon a^{\frac{p-1+\theta}{p-1}}+C_{\epsilon,p,\theta}b^{\frac{p-1+\theta}{\theta}}$ for any $\epsilon>0$, we get for $t\leq T_p$,
\begin{align}\label{4.24}
 & X_n(t) + \frac{p}{2}E\int_0^{t\wedge\tau_M^n} \Vert u_n(s)\Vert^{p-2}\Vert u_n(s)\Vert_{V}^2ds \nonumber\\
 \leq & e^{C (e^{T_p} -1)} E\left\{\left[M(t\wedge\tau_M^n) +1\right]^{e^{T_p}}\right\} = C E\left\{\left[M(t\wedge\tau_M^n) +1\right]^{\frac{p}{p-1+\theta}}\right\} \nonumber\\
  \leq & C\left(1+ \Vert u_0\Vert^{\frac{p^2}{p-1+\theta}}\right)  \nonumber\\
  & + CE\sup_{r\in[0,t\wedge\tau_M^n]}\left|\int_0^r \Vert u_n(s)\Vert^{p-2}\big(\sigma(u_n(s)),u_n(s)\big) dB_s \right|^{\frac{p}{p-1+\theta}} \nonumber\\
  \leq & C\left(1+ \Vert u_0\Vert^{\frac{p^2}{p-1+\theta}}\right)+  CE\left(\int_0^{t\wedge\tau_M^n}\Vert u_n(s)\Vert^{2p-2}\Vert\sigma(u_n(s))\Vert^2  ds\right)^{\frac{p}{2(p-1+\theta)}} \nonumber\\
  \leq & C\left(1+ \Vert u_0\Vert^{\frac{p^2}{p-1+\theta}}\right) \nonumber\\
  & +CE\left[\sup_{s\in[0,t\wedge\tau_M^n]}\Vert u_n(s)\Vert^{(2p-2)\times\frac{p}{2(p-1+\theta)}}\times \left(\int_0^{t\wedge\tau_M^n}\Vert u_n(s)\Vert^{2\theta} ds\right)^{\frac{p}{2(p-1+\theta)}}\right] \nonumber\\
  \leq &  C\left(1+ \Vert u_0\Vert^{\frac{p^2}{p-1+\theta}}\right)+ \epsilon X_n(t) + C_{\epsilon} E\left(\int_0^{t\wedge\tau_M^n}\Vert u_n(s)\Vert^{2\theta} ds\right)^{\frac{p}{2\theta}} \nonumber\\
  \leq &  C\left(1+ \Vert u_0\Vert^{\frac{p^2}{p-1+\theta}}\right)+ \epsilon X_n(t) + C_{\epsilon} E\int_0^{t\wedge\tau_M^n} \Vert u_n(s)\Vert^{p} ds  \nonumber\\
  \leq &  C\left(1+ \Vert u_0\Vert^{\frac{p^2}{p-1+\theta}}\right)+ \epsilon X_n(t) + C_{\epsilon}\int_0^t X_n(s)ds .
  \end{align}
Note that the above constants $C, C_{\epsilon}$ are independent of $M, n$. Hence subtracting $\epsilon X_n(t)$ from both sides of above inequality, then applying Gronwall's inequality, we obtain for any $n\in\mathbb{N}$ and $M>0$,
\begin{align}\label{4.25}
  X_n (T_p)+ E\int_0^{T_p \wedge\tau_M^n} \Vert u_n(s)\Vert^{p-2}\Vert u_n(s)\Vert_{V}^2ds  \leq C_{p,\theta}\left(1+ \Vert u_0\Vert^{\frac{p^2}{p-1+\theta}}\right).
\end{align}
Therefore, letting $M\rightarrow\infty$ and by Fatou's lemma, (\ref{4.16}) is proved. $\blacksquare$

\vskip 0.5cm

Consider the time interval $[0,T_2]$ (see (\ref{Tp})).
Take $\beta\in(0,1)$ and $p>1$. Let $W^{\beta, p}([0,T_2];V^*)$ be the space of functionals $u(\cdot): [0,T_2]\rightarrow V^*$ with the finite norm defined by
\begin{align}\label{4.27}
  \Vert u\Vert^{p}_{{W^{\beta, p}([0,T_2];V^*)}}:= \int_0^{T_2} \Vert u(t)\Vert^p_{V^*} dt + \int_0^{T_2}\int_0^{T_2} \frac{\Vert u(t)-u(s)\Vert_{V^*}^p}{|t-s|^{1+\beta p}} dtds .
\end{align}

\begin{lemma}\label{190630.0140}
Suppose hypothesis (\hyperlink{H.2}{H.2}) holds.
If $\beta<\frac{1}{2}$ and $1<p<2$, then
\begin{align}\label{4.28}
  \sup_{n} E\left(\Vert u_n\Vert_{W^{\beta, p}([0,T_2];V^*)}^p\right)< \infty .
\end{align}
\end{lemma}
\noindent {\bf Proof}.
 Consider
\begin{align}\label{4.29}
  u_n(t)-u_n(0)& =\int_0^t \Delta u_n(r) dr + \int_0^t u_n(r)\log|u_n(r)| dr + \int_0^t \sigma(u_n(r))dB_r \nonumber\\
  & =: J_1(t) +J_2(t) +J_3(t) .
\end{align}
We have
\begin{align}\label{4.30}
 E\Vert u_n(t)-u_n(s)\Vert_{V^*}^p
  \leq & 3^{p-1}\times(\Vert J_1(t)-J_1(s)\Vert_{V^*}^p  +\Vert J_2(t)-J_2(s)\Vert_{V^*}^p \nonumber\\
  & ~~~~~~~~~+\Vert J_3(t)-J_3(s)\Vert_{V^*}^p) .
\end{align}
Without loss of generality, we can assume that $t\geq s$. By H\"o{}lder's inequality, we have
\begin{align}\label{4.31}
  E\Vert J_1(t)-J_1(s)\Vert_{V^*}^p\leq & CE\left(\int_s^t \Vert u_n(r)\Vert_{V} dr\right)^p \nonumber\\
  \leq & C E\left(1+\int_0^{T_2} \Vert u_n(r)\Vert_{V}^2 dr\right) \times |t-s|^{\frac{p}{2}} .
\end{align}
From $H^1\hookrightarrow L^q$, it follows that $L^{q^*}\hookrightarrow V^*$, where $q*\in[\frac{2d}{d+2},2)$ if the dimension $d>2$, and $q^*\in(1,2)$ if $d=1,2$.
Note that for any $\epsilon>0$, there exists a constant $C_{\epsilon}>0$ such that for any $a\geq 0$,
\begin{align*}
  \big|a\log|a|\big|\leq C_{\epsilon}\left(1+a^{1+\epsilon}\right) .
\end{align*}
Now we take $\epsilon>0$ sufficiently small so that $(1+\epsilon)p\leq 2$ and $(1+\epsilon)q^*\leq 2$.
Thus
\begin{align}\label{4.32}
  E\Vert J_2(t)-J_2(s)\Vert_{V^*}^p\leq & CE\left(\int_s^t \Vert u_n(r)\log|u_n(r)|\Vert_{L^{q^*}} dr\right)^p \nonumber\\
  \leq & C_{\epsilon}E\left[\int_s^t \big(1+\Vert u_n(r)\Vert^{1+\epsilon}\big) dr\right]^p \nonumber\\
  \leq & C\left(1+E\sup_{r\in[0,T_2]}\Vert u_n(r)\Vert^2\right)\times |t-s|^p .
\end{align}
Similarly,
\begin{align}\label{4.33}
  E\Vert J_3(t)-J_3(s)\Vert_{V^*}^p\leq & CE\left(\int_s^t \Vert \sigma(u_n(r))\Vert_{V^*}^2 dr\right)^{\frac{p}{2}} \nonumber\\
  \leq & C\left(1+E\sup_{r\in[0,T_2]}\Vert u_n(r)\Vert^2\right)\times |t-s|^{\frac{p}{2}} .
\end{align}
Putting (\ref{4.30})-(\ref{4.33}) together, and in view of (\ref{4.27}), we see that
\begin{align}\label{4.34}
   E\Vert u_n\Vert^{p}_{{W^{\beta, p}([0,T_2];V^*)}}
  \leq &  C\left(1+E\sup_{r\in[0,T_2]}\Vert u_n(r)\Vert^2 + E\int_0^{T_2} \Vert u_n(r)\Vert_{V}^2 dr\right) \nonumber\\
  & \times \left(1+\int_0^{T_2}\int_0^{T_2} \frac{|t-s|^{\frac{p}{2}}}{|t-s|^{1+\beta p}} dtds \right) .
\end{align}
Therefore, if $\beta<\frac{1}{2}$ then the above integral is finite. Moreover, (\ref{4.28}) follows from (\ref{4.16}).
$\blacksquare$

\begin{lemma}\label{lemma 4.5}
Suppose hypothesis (\hyperlink{H.2}{H.2}) holds. Then for any $1<p<2$,  $\{u_n\}$ is tight in $L^p([0,T_2];H)$ and in $C([0,T_2];V^*)$.
\end{lemma}
\noindent {\bf Proof}.
It is known that $L^p([0,T_2];V)\cap W^{\beta, p}([0,T_2];V^*)$ is compactly embedded into $L^p([0,T_2];H)$, see Theorem 2.1 in \cite{FG}. So for any $L>0$, the set
\begin{align*}
  K_L :=\left\{u\in L^p([0,T_2];H): \Vert u\Vert_{L^p([0,T_2];V)}+\Vert u \Vert_{W^{\beta, p}([0,T_2];V^*)}\leq L\right\}
\end{align*}
is a compact subset of $L^p([0,T_2];H)$. Moreover,
\begin{align}\label{4.36}
  & \lim_{L\rightarrow\infty}\sup_{n}P(u_n\notin K_L) \nonumber\\
  =& \lim_{L\rightarrow\infty}\sup_{n}P\left(\Vert u\Vert_{L^p([0,T_2];V)}+\Vert u \Vert_{W^{\beta, p}([0,T_2];V^*)}> L\right) \nonumber\\
  \leq & \lim_{L\rightarrow\infty}\frac{2^{p-1}}{L^p}\sup_n E\left(\Vert u_n\Vert_{L^p([0,T_2];V)}^p+\Vert u_n \Vert_{W^{\beta, p}([0,T_2];V^*)}^p\right) \nonumber\\
  = & 0 .
\end{align}
Hence $\{u_n\}$ is tight in $L^p([0,T_2];H)$.

\vskip 0.3cm

Now we prove the tightness in $C([0,T_2];V^*)$. Note that $H$ is compactly embedded into $V^*$.
One the one hand,
\begin{align}\label{4.37}
  & \lim_{L\rightarrow\infty}\sup_{n}P\left(\sup_{t\in[0,T_2]}\Vert u_n(s)\Vert >L \right) \nonumber\\
  \leq & \lim_{L\rightarrow\infty} \frac{1}{L^2}\sup_{n}E\left[\sup_{t\in[0,T_2]}\Vert u_n(t)\Vert^2 \right]=0 .
\end{align}
One the other hand,
for any stopping time $0\leq\zeta_{n}\leq T_2$ and any $\varepsilon>0$, Similar to (\ref{4.30})-(\ref{4.33}), we have
\begin{align}\label{4.38}
  & \lim_{\delta\rightarrow 0}\sup_{n}P(\Vert u_{n}(\zeta_{n}+\delta) -u_{n}(\zeta_{n})\Vert_{V^*}>\varepsilon) \nonumber\\
  \leq & \frac{1}{\varepsilon^p}\lim_{\delta\rightarrow 0}\sup_{n}E\Vert u_{n}(\zeta_{n}+\delta) -u_{n}(\zeta_{n})\Vert_{V^*}^p \nonumber\\
  \leq & \frac{C}{\varepsilon^p}\lim_{\delta\rightarrow 0}\sup_{n}\left[\left(1+E\sup_{r\in[0,T_2]}\Vert u_n(r)\Vert^2 + E\int_0^{T_2} \Vert u_n(r)\Vert_{V}^2 dr\right)\times\delta^{\frac{p}{2}}\right] \nonumber\\
  = & 0 ,
\end{align}
where $\zeta^{\varepsilon}+\delta :=T_2 \wedge(\zeta^{\varepsilon}+\delta)$. By Aldou's tightness criterion (see Theorem 1 in \cite{A}), $\{u_n\}$ is tight in $C([0,T_2];V^*)$. $\blacksquare$

%
%
%
%
%

\subsection{Existence of solutions}
Here is the main result in this section.

\begin{theorem}\label{theorem 4.6}
Suppose hypotheses (\hyperlink{H.1}{H.1}) and (\hyperlink{H.2}{H.2}) hold. Then there exists a unique global solution $u$ to (\ref{Abstract}) for every initial value $u_0\in L^2(D)$. Moreover, for any $T>0$ and $p\geq 2$, we have
\begin{align}\label{4.39}
     E\sup_{t\in[0,T]} \Vert u(t)\Vert^p + E\int_0^{T} \Vert u\Vert^{p-2}\Vert u\Vert_{V}^2 ds
  <\infty .
\end{align}
\end{theorem}
\noindent {\bf Proof}.
The proof is divided into two steps. We first prove the existence of the solutions in the time interval $[0,T_2]$. Then we piece together solutions on subintervals to get  the global existence of the solution and also the estimate (\ref{4.39}).

\vskip 0.4cm

Step 1. We will prove the existence of a probabilistic weak solution in the time interval $[0,T_2]$ and then use the Yamada-Watanabe theorem to obtain the strong solution.
\vskip 0.3cm

Fix  $1<r<2$, and set
\[\Upsilon:= [ L^r ([0,T_2];H)\cap C([0,T_2];V^*) ] \times C([0,T_2];\mathbb{R}) .\]
From Lemma \ref{lemma 4.5}, we see that the law $\mathcal{L}(u_n, B)$ of the random vector $(u_n, B)$ is tight in $\Upsilon$. By Prokhorov's theorem, there exist a subsequence, still denoted by  $(u_n , B)$, and a probability measure $\mu$ on $\Upsilon$, such that $\mathcal{L}(u_{n}, B)$ is weakly convergent to $\mu$.

By the modified version of the Skorokhod embedding theorem whose proof can be found in Appendix C of \cite{BHR}, there exist a new probability space $(\widetilde{\Omega}, \widetilde{\mathcal{F}},\widetilde{P})$ and a sequence of $\Upsilon$-valued random vectors $\{(\widetilde{u}_n, \widetilde{B}_n)\}$ and $(\widetilde{u}, \widetilde{B})$ such that for any $n\in\mathbb{N}$, $\widetilde{B}_n=\widetilde{B}$, $\widetilde{P}$-a.s., and $\mathcal{L}(\widetilde{u}_n, \widetilde{B}_n)=\mathcal{L}(u_n, B)$, $\mathcal{L}(\widetilde{u}, \widetilde{B})=\mu$, moreover, $\widetilde{P}$-a.s.,
\begin{align}\label{exist 1}
  \Vert \widetilde{u}_n -\widetilde{u}\Vert_{L^r ([0,T_2];H)}+ \Vert \widetilde{u}_n -\widetilde{u}\Vert_{C([0,T_2];V^*)}\rightarrow 0 .
\end{align}

From the equation satisfied by the random vector $(u_n, B)$, we see that $(\widetilde{u}_n, \widetilde{B})$ satisfies the following equation in $V^*$:
\begin{align}\label{0403.6}
  \widetilde{u}_n(t)=& P_n u_0 + \int_0^t \Delta \widetilde{u}_n(s) ds +\int_0^t P_n[\widetilde{u}_n(s)\log|\widetilde{u}_n(s)| ]ds  \nonumber\\
  & +\int_0^t P_n\sigma(\widetilde{u}_n(s))d\widetilde{B}_s,\quad t\in[0,T_2] ,
\end{align}
where $P_n$ is the projection operator defined in (\ref{4.1}).
Therefore, $\{\widetilde{u}_n\}$ also satisfies
\begin{align}\label{0403.1}
    \sup_{n}\widetilde{E}\left[\sup_{s\in[0,T_2]}\Vert \widetilde{u}_n(s)\Vert^2 + \int_0^{T_2} \Vert \widetilde{u}_n(s)\Vert_{V}^2 ds \right] <\infty .
\end{align}
From $\Vert \widetilde{u}_n -\widetilde{u}\Vert_{C([0,T_2];V^*)}\rightarrow 0$, $\widetilde{P}$-a.s., it follows that
\begin{align}\label{0403.2}
  \Vert P_m\widetilde{u}_n-P_m \widetilde{u}\Vert_{C([0,T_2];V)}\xrightarrow{n\rightarrow \infty} 0 , \quad \forall\, m\in\mathbb{N} .
\end{align}
Hence by Fatou's lemma and (\ref{0403.1}), we have
\begin{align}\label{0403.4}
  \widetilde{E}\sup_{t\in[0,T_2]}\Vert  \widetilde{u}(t)\Vert^2 \leq & \widetilde{E}\liminf_{m\rightarrow\infty}\sup_{t\in[0,T_2]}\Vert P_m \widetilde{u}(t)\Vert^2 \nonumber\\
  \leq & \liminf_{m\rightarrow\infty}\widetilde{E}\sup_{t\in[0,T_2]}\Vert P_m \widetilde{u}(t)\Vert^2 \nonumber\\
  \leq & \liminf_{m\rightarrow\infty}\liminf_{n\rightarrow\infty}\widetilde{E}\sup_{t\in[0,T_2]}\Vert P_m \widetilde{u}_n(t)\Vert^2 \nonumber\\
  \leq & \sup_{n}\widetilde{E}\sup_{t\in[0,T_2]}\Vert \widetilde{u}_n(t)\Vert^2 <\infty .
\end{align}
Similarly,
\begin{align}\label{0403.5}
  \widetilde{E}\int_0^{T_2} \Vert \widetilde{u}(s)\Vert_{V}^2 ds <\infty .
\end{align}
By (\ref{exist 1}), (\ref{0403.1}) and the sublinear growth of $\sigma$, we see that
there exists a subsequence, still denoted by  $\{\widetilde{u}_{n}\}$, such that as $n\rightarrow\infty$,
\begin{itemize}
  \item [(i)] $\widetilde{u}_n(\widetilde{\omega},t,x)\rightarrow \widetilde{u}(\widetilde{\omega},t,x)$, a.e. $(\widetilde{\omega},t,x)\in\widetilde{\Omega}\times [0,T_2]\times D$;
  \item [(ii)]  $\widetilde{u}_n\rightarrow \widetilde{u}$ strongly in $L^r(\Omega;L^r([0,T_2];H))$;
  \item [(iii)] $\widetilde{u}_n\rightharpoonup \widetilde{u}$ weakly in $L^2(\widetilde{\Omega};L^2([0,T_2];V))$;
  \item [(iv)] $\Delta\widetilde{u}_n \rightharpoonup\Delta\widetilde{u}$ weakly in $L^2(\widetilde{\Omega};L^2([0,T_2];V^*))$;
  \item [(v)] $P_n[\widetilde{u}_n\log|\widetilde{u}_n|] \rightarrow \widetilde{u}\log|\widetilde{u}|$ strongly in $L^r(\widetilde{\Omega};L^r([0,T_2];V^*))$;
  \item [(vi)] $\int_0^{\cdot} P_n\sigma(\widetilde{u}_n(s))d\widetilde{B}_s\rightarrow\int_0^{\cdot} \sigma(\widetilde{u}(s))d\widetilde{B}_s$ strongly in $L^{\infty}([0,T_2];L^2(\widetilde{\Omega};V^*))$.
\end{itemize}
Letting $n\rightarrow\infty$ in (\ref{0403.6}), it is easy to see that $\widetilde{u}$ satisfies
\begin{align}\label{0403.7}
  \widetilde{u}(t)=u_0 + \int_0^t \Delta \widetilde{u}(s) ds +\int_0^t \widetilde{u}(s)\log|\widetilde{u}(s)| ds  +\int_0^t \sigma(\widetilde{u}(s))d\widetilde{B}_s
\end{align}
$\widetilde{P}$-a.s., for any $t\in[0,T_2]$.
Note that $\Vert\widetilde{u}(s)\log|\widetilde{u}(s)|\Vert\leq C(1+\Vert \widetilde{u}(s)\Vert_{V}^2)$. The continuity of $\widetilde{u}$ as an 
$H$-valued process follows from the above equation.
Hence we obtain the existence of a probabilistic weak solutions on the time interval $[0,T_2]$.

The existence of probabilistic strong solutions on $[0,T_2]$ follows from the pathwise uniqueness of the  solutions proved in Section 3 and the Yamada-Watanabe theorem.


\vskip 0.5cm

Step 2. We prove the global existence of the probabilistic strong solutions and the moment estimates.

\vskip 0.3cm

In the following, we take any fixed $T>0$ and any fixed $p\geq 2$. We will construct a solution on the time interval $[0,T]$, and establish the estimate (\ref{4.39}).
First, we introduce some notations. For $z\in [2,\infty)$, define functions
\begin{gather}
  \gamma^0(z):=z,\quad \gamma(z):=\frac{z^2}{z-1+\theta},\quad \gamma^2(z):= \gamma(\gamma(z)), \nonumber\\
  \gamma^n(z):=\underbrace{\gamma(\gamma(\dots\gamma(z)\dots))}_{\text{nth iteration}},\quad n\in\mathbb{N} .
\end{gather}
For any $z\geq 2$, we have $\gamma(z)-z>1-\theta $. So for any $z\geq 2$, $\gamma^n(z)\uparrow\infty$ as $n\rightarrow\infty$.
If $z\geq \frac{1}{\theta}-1$, then $\gamma(z)\leq z+1$.
Let
\begin{align}
  i_p:=\min \left\{i\geq 0: \gamma^i(p)\geq \frac{1}{\theta}-1 \right\} .
\end{align}
To simplify the notation, we denote
\begin{align}
  q(i):=\gamma^{i}(p), \quad T(i):=T_{q(i)}=T_{\gamma^{i}(p)},\quad i\in\mathbb{N}.
\end{align}
Note that for $i\geq i_p$, $ \gamma^{i}(p)\leq \gamma^{i_p}(p)+i-i_p$, because $\gamma(z)\leq z+1$.
Then we have
\begin{align}\label{4.49}
  \sum_{i=0}^{\infty} T(i)\geq & \sum_{i=i_p +1}^{\infty} T_{\gamma^i(p)}
   \geq \sum_{i=i_p +1}^{\infty} T_{\gamma^{i_p}(p)+i-i_p} \nonumber\\
   = & \sum_{j=1}^{\infty} \log \frac{\gamma^{i_p}(p)+j}{\gamma^{i_p}(p)+j -1 +\theta} \nonumber\\
   = &  \sum_{j=1}^{\infty} \log \left( 1+\frac{1-\theta}{\gamma^{i_p}(p)+j -1 +\theta}\right) \nonumber\\\
   \geq &\delta \sum_{j=1}^{\infty} \frac{1-\theta}{\gamma^{i_p}(p)+j -1 +\theta} =\infty ,
\end{align}
where $\delta$ is a small positive constant.
%
%
%
Define
\begin{align}
  \kappa:=\min\left\{n\geq 0: \sum_{i=0}^n T(i)\geq T\right\} .
\end{align}
From (\ref{4.49}), we see that $\kappa<\infty$.  Set $S(0)=0$, and
\begin{align}
  S(i+1)=S(i)+T(\kappa-i),\quad i=0,1,\dots ,\kappa.
\end{align}
Then by the definition of $\kappa$,
\[S(\kappa)<T\leq S(\kappa+1) .\]

For any $S\geq 0$, define
\begin{align}
  B^{S}_t:=B(t+S)-B(S) ,\quad \mathcal{F}^{S}_t:=\mathcal{F}_{t+S}, \quad t\geq 0 .
\end{align}
Then $(B^{S}_t)_{t\geq 0}$ is a Brownian motion with respect to the filtration $(\mathcal{F}^{S}_t)_{t\geq 0}$.
Consider the equation
\begin{align}
  u^{S,h}(t)= & h+\int_0^t \Delta u^{S,h}(s) ds +\int_0^t u^{S,h}(s)\log|u^{S,h}(s)| ds \nonumber\\
  & +\int_0^t \sigma(u^{S,h}(s))dB^{S}_s , \quad t\in [0,T_2] .
\end{align}
Then by Step 1, we know that there exists a unique probabilistic strong solution $\{u^{S,h}(t)\}_{t\in[0,T_2]}$ to the above equation. Furthermore, $u^{S,h}$ is a measurable map of $h\in H$.



Now, by iteration we define a process
\begin{align}
      u(t):=u^{S(i),u(S(i))}(t-S(i)),\quad t\in[S(i),S(i+1)],\quad i=0,1,\dots ,\kappa.
\end{align}
Since $u(S(i))$ is an $\mathcal{F}_{S(i)}$ measurable random variable, which is independent of the Brownian motion $B^{S(i)}$,
we see that $u$ satisfies the following equation
\begin{align}
  u(t)= & u(S(i)) + \int_{S(i)}^t \Delta u(s) ds +\int_{S(i)}^t u(s)\log|u(s)| ds \nonumber\\
  & +\int_{S(i)}^t \sigma(u(s))dB_s , \quad t\in[S(i),S(i+1)] .
\end{align}
This means that the  process $u$ satisfies equation (\ref{Abstract}) in the time interval $[0,S(\kappa+1)]$. Furthermore, $u$ is an $H$-valued continuous and $(\mathcal{F}_t)$-adapted process and
\[
\int_0^{S(\kappa+1)} \Vert u(s)\Vert_{V}^2 ds <\infty , \quad P\text{-a.s.}
\]
Note that $T\leq S(\kappa+1)$.
Therefore, for any $T> 0$, we have constructed a probabilistic strong solution to equation (\ref{Abstract}) on the interval $[0,T]$. By the arbitrariness of $T$ and the pathwise uniqueness, the solution $u$ is global.

\vskip 0.3cm

Next we estimate the moment of $u$.
Since we have shown that the estimates (\ref{0403.4}) and (\ref{0403.5}) are satisfied for the probabilistic weak solution on $[0,T_2]$, they are also satisfied for probabilistic strong solutions on $[0,T_2]$. Similarly, we have for any $i,j\in\{0,\dots, \kappa\}$ and $j\leq i$,
\begin{align}
  & E\sup_{t\in[0,T(i)]}\Vert u^{0,h}(t)\Vert^{q(j)} + E\int_0^{T(i)} \Vert u^{0,h}(t)\Vert^{q(j)-2}\Vert u^{0,h}(t)\Vert_V^2 dt\nonumber\\
  \leq & C_{p,\theta} \left(1+\Vert h\Vert^{q(j+1)}\right) .
\end{align}
In particular,
\begin{align}
    \Phi_i(h):=& E\sup_{t\in[0,T(i)]}\Vert u^{0,h}(t)\Vert^{q(i)} + E\int_0^{T(i)} \Vert u^{0,h}(t)\Vert^{p-2}\Vert u^{0,h}(t)\Vert_V^2dt \nonumber\\
    \leq & C_{p,\theta} \left(1+\Vert h\Vert^{q(i+1)}\right) .
\end{align}
Using the fact that $u(S(i))$ is independent of the Brownian motion $(B^{S(i)}_t)_{t\geq 0}$, and the process $u^{S,h}$ has the same distribution with the process $u^{0,h}$, we get
\begin{align}
  & E\sup_{t\in[S(i),S(i+1)]} \Vert u(t)\Vert^{q(\kappa-i)}+E\int_{S(i)}^{S(i+1)}\Vert u(t)\Vert^{p-2}\Vert u(t)\Vert_{V}^2 dt \nonumber\\
  = & E \sup_{t\in[0,T(\kappa-i)]} \Vert u^{S(i),u(S(i))}(t)\Vert^{q(\kappa-i)} \nonumber\\
  & +E\int_{0}^{T(\kappa-i)}\Vert u^{S(i),u(S(i))}(t)\Vert^{p-2}\Vert u^{S(i),u(S(i))}(t) \Vert_{V}^2 dt \nonumber\\
  = & E\Bigg[\Bigg( E \sup_{t\in[0,T(\kappa-i)]} \Vert u^{S(i),h}(t)\Vert^{q(\kappa-i)} \nonumber\\
     & +E\int_{0}^{T(\kappa-i)}\Vert u^{S(i),h}(t)\Vert^{p-2}\Vert u^{S(i),h}(t) \Vert_{V}^2 dt\Bigg)\Bigg|_{h=u(S(i))}\Bigg] \nonumber\\
  = & E\Phi_{\kappa-i}(u(S(i))) \nonumber\\
  \leq & C_{p,\theta}(1+E\Vert u(S(i))\Vert^{q(\kappa-i+1)}) .
\end{align}
Therefore, by mathematical induction, we obtain
\begin{align}
  & E\sup_{t\in[0,S(\kappa+1)]} \Vert u(t)\Vert^p + E\int_0^{S(\kappa+1)} \Vert u\Vert^{p-2}\Vert u\Vert_{V}^2 ds \nonumber\\
  \leq & C_{p,\theta} (1+\Vert u_0\Vert^{q(\kappa+1)}) .
\end{align}
Note that $S(\kappa)<T\leq S(\kappa+1)$, so the estimate (\ref{4.39}) is deduced.
The proof of Theorem \ref{theorem 4.6} is complete.  $\blacksquare$

%
%

%

\section{Existence of solutions: part II}
\setcounter{equation}{0}
 In this section, we assume (\hyperlink{H.3}{H.3}) holds.  We will modify the arguments in Section $5$ to get the existence of the solutions.
\vskip 0.5cm
First
we note that under hypothesis (\hyperlink{H.3}{H.3}), there still exists a global solution $u_n$ to the approximating equation (\ref{4.2}) according to Theorem \ref{190630.0008}. To establish the tightness of $\{u_n, n\geq 1\}$ under this new condition, we need some new estimates. Define
\begin{align*}
\rho(x) :=
\begin{cases}
\log x  & x \in[\mathrm{e},\infty) , \\
\frac{x}{\mathrm{e}} &  x\in[0, \mathrm{e}] .
\end{cases}
\end{align*}
and
\[ \Phi(z):=\exp\left (\int_0^z\frac{1}{1+x + x\rho(x)}dx\right ), \quad  z\in \mathbb{R}_{+}. \]
Then
\begin{align}\label{190630.0118.1}
  \Phi^{\prime}(z)=\Phi(z)\times\frac{1}{1+z+ z\rho(z)} ,  \quad \Phi^{\prime\prime}(z)\leq 0 .
\end{align}
For any $M>0$, introduce stopping times
\begin{equation}\label{5.1}
  \tau_M^n:= \inf\{t>0: \Vert u_n(t)\Vert^2>M\}.
\end{equation}
Set $u_n^M(t)=u_n(t\wedge \tau_M^n)$.
The following lemma is crucial.
\begin{lemma}\label{190630.0143}
Suppose hypothesis (\hyperlink{H.3}{H.3}) holds. Then there is a constant $C$ such that for all $t, M>0$,
\begin{equation}\label{5.2}
  \sup_n E\left[\Phi(\Vert u_n^M(t)\Vert^2) + \int_0^{t\wedge\tau_M^n}\Phi^{\prime}(\Vert u_n^M(s)\Vert^2) \Vert u_n^M(s)\Vert_{V}^2 ds\right]\leq \Phi(\Vert u_0\Vert^2)e^{Ct}.
\end{equation}
\end{lemma}
\noindent{\bf Proof}. The proof is inspired by the proof of Theorem A in \cite{FZ}.  As in the proof Lemma \ref{lemma 190706.1431}, we have
\begin{align}\label{5.3}
  \Vert u_n^M(t)\Vert^2 = & \Vert u_0 \Vert^2 - 2\int_0^{t\wedge\tau_M^n} \Vert u_n^M(s)\Vert_{V}^2 ds \nonumber\\
  & + 2\int_0^{t\wedge\tau_M^n}  \big(u_n^M(s)\log|u_n^M(s)|,u_n^M(s)\big) ds \nonumber\\
  & +  \int_0^{t\wedge\tau_M^n} \Vert P_n\sigma(u_n^M(s))\Vert^2 ds \nonumber\\
  & + 2\int_0^{t\wedge\tau_M^n} \big(\sigma(u_n^M(s)),u_n^M(s)\big)dB_s.
\end{align}
Applying Ito's formula to the real valued process $\Vert u_n^M\Vert^2$, we obtain
\begin{align}\label{5.4}
  \Phi(\Vert u_n^M(t)\Vert^2) = &\Phi(\Vert u_0 \Vert^2) - 2\int_0^{t\wedge\tau_M^n}\Phi^{\prime}(\Vert u_n^M(s)\Vert^2) \Vert u_n^M(s)\Vert_{V}^2 ds \nonumber\\
  & + 2\int_0^{t\wedge\tau_M^n} \Phi^{\prime}(\Vert u_n^M(s)\Vert^2) \big(u_n^M(s)\log|u_n^M(s)|,u_n^M(s)\big) ds \nonumber\\
  & +  \int_0^{t\wedge\tau_M^n} \Phi^{\prime}(\Vert u_n^M(s)\Vert^2)\Vert P_n\sigma(u_n^M(s))\Vert^2 ds \nonumber\\
  & + 2\int_0^{t\wedge\tau_M^n} \Phi^{\prime}(\Vert u_n^M(s)\Vert^2)\big(\sigma(u_n^M(s)),u_n^M(s)\big)dB_s\nonumber\\
  &+2\int_0^{t\wedge\tau_M^n} \Phi^{\prime\prime}(\Vert u_n^M(s)\Vert^2)\big(\sigma(u_n^M(s)),u_n^M(s)\big)^2ds.
\end{align}
From the logarithmic Soblev inequality (\ref{log sobolev}) with $\varepsilon=\frac{1}{4}$ and hypothesis (\hyperlink{H.3}{H.3}), it follows that
\begin{align}\label{5.5}
  & \Phi(\Vert u_n^M(t)\Vert^2) \nonumber\\
  \leq &\Phi(\Vert u_0 \Vert^2) - 2\int_0^{t\wedge\tau_M^n}\Phi^{\prime}(\Vert u_n^M(s)\Vert^2) \Vert u_n^M(s)\Vert_{V}^2 ds \nonumber\\
  & + 2\int_0^{t\wedge\tau_M^n} \Phi^{\prime}(\Vert u_n^M(s)\Vert^2) \bigg(\frac{1}{4}\Vert u_n^M(s)\Vert_{V}^2+\frac{d\log 4}{4}\Vert u_n^M(s)\Vert^2 \nonumber \\
  & ~~~~~~~~~~~~~~~~~~~~~~~~~~~~~~~~~~+ \Vert u_n^M(s)\Vert^2\log\Vert u_n^M(s)\Vert\bigg) ds\nonumber\\
  & +  \int_0^{t\wedge\tau_M^n} \Phi^{\prime}(\Vert u_n^M(s)\Vert^2) \int_D \left(2C_3^2 + 2C_4^2 |u_n^M(s,x)|^2\log_{+}|u_n^M(s,x)| \right)dx ds \nonumber\\
  & + 2\int_0^{t\wedge\tau_M^n} \Phi^{\prime}(\Vert u_n^M(s)\Vert^2)\big(\sigma(u_n^M(s)),u_n^M(s)\big)dB_s\nonumber\\
  &+2\int_0^{t\wedge\tau_M^n} \Phi^{\prime\prime}(\Vert u_n^M(s)\Vert^2)\big(\sigma(u_n^M(s)),u_n^M(s)\big)^2ds.
\end{align}
Applying the modified logarithmic Soblev inequality (\ref{log sobolev modi}) with $\varepsilon=\frac{1}{4C_4^2}$, we deduce from (\ref{5.5}) that
\begin{align}\label{5.6}
&\Phi(\Vert u_n^M(t)\Vert^2)+\int_0^{t\wedge\tau_M^n}\Phi^{\prime}(\Vert u_n^M(s)\Vert^2) \Vert u_n^M(s)\Vert_{V}^2 ds\nonumber\\
  \leq & \Phi(\Vert u_0 \Vert^2) \nonumber\\
   & + C\int_0^{t\wedge\tau_M^n} \Phi^{\prime}(\Vert u_n^M(s)\Vert^2) \Big(1+ \Vert u_n^M(s)\Vert^2+ \Vert u_n^M(s)\Vert^2\log\Vert u_n^M(s)\Vert\Big) ds\nonumber\\
  & + 2\int_0^{t\wedge\tau_M^n} \Phi^{\prime}(\Vert u_n^M(s)\Vert^2)\big(\sigma(u_n^M(s)),u_n^M(s)\big)dB_s\nonumber\\
  &+2\int_0^{t\wedge\tau_M^n} \Phi^{\prime\prime}(\Vert u_n^M(s)\Vert^2)\big(\sigma(u_n^M(s)),u_n^M(s)\big)^2 ds\nonumber\\
  =& : \Phi(\Vert u_0 \Vert^2)+I_1^n(t)+I_2^n(t)+I_3^n(t).
  \end{align}
  We now estimate each of the terms on the right side of the above inequality. From (\ref{190630.0118.1}) we see that
  \begin{align}\label{5.7}
  I_1^n(t)= & C\int_0^{t\wedge\tau_M^n} \Phi(\Vert u_n^M(s)\Vert^2) \frac{1}{1+\Vert u_n^M(s)\Vert^2+ \Vert u_n^M(s)\Vert^2 \rho(\Vert u_n^M(s)\Vert^2)} \nonumber\\
  &~~~~~~\times  \Big(1+\Vert u_n^M(s)\Vert^2+ \Vert u_n^M(s)\Vert^2\log\Vert u_n^M(s)\Vert\Big) ds\nonumber\\
  \leq & C\int_0^t \Phi(\Vert u_n^M(s)\Vert^2)ds.
  \end{align}
  Note that $\Phi^{\prime\prime}(z)\leq 0$  for all $z\geq 0$. Therefore, $I_3^n(t)\leq 0$. Now substitute (\ref{5.7})
  into (\ref{5.6}) and take expectation on both sides of (\ref{5.6}) to obtain
  \begin{align}\label{5.8}
  & E[\Phi(\Vert u_n^M(t)\Vert^2)] + E\int_0^{t\wedge\tau_M^n}\Phi^{\prime}(\Vert u_n^M(s)\Vert^2) \Vert u_n^M(s)\Vert_{V}^2 ds \nonumber\\
  \leq & \Phi(\Vert u_0 \Vert^2)+ C\int_0^t E[\Phi(\Vert u_n^M(s)\Vert^2)]ds.
  \end{align}
An application of  the Gronwall inequality yields the desired result. $\blacksquare$
\begin{lemma}\label{190630.0142}
Suppose hypothesis (\hyperlink{H.3}{H.3}) holds. Then for any $T, M>0$,
\begin{align}\label{5.13}
  \sup_n\sup_{0\leq t\leq T}\Vert u_n^M(t)\Vert^2 \leq M , \quad P\text{-a.s.},
\end{align}
moreover, there exists a constant $C_{T, M}$ such that
\begin{equation}\label{5.9}
\sup_n E\left[\int_0^{T\wedge\tau_M^n}\Vert u_n^M(t)\Vert^2_Vdt\right ]\leq C_{T,M} .
\end{equation}
\end{lemma}
\noindent{\bf Proof}. According to the definition of $\tau_M^n$, (\ref{5.13}) is obvious.
By (\ref{190630.0118.1}), it can be seen that
\begin{align}
 c_M:=\inf_{0\leq z\leq M} \Phi^{\prime}(z)>0
\end{align}
Hence from (\ref{5.2}) and (\ref{5.13}) it follows that
\begin{align}
  c_M \times\sup_n E\int_{0}^{T\wedge\tau_M^n} \Vert u_n^M(s)\Vert_{V}^2 ds \leq & \sup_n E \int_0^{T\wedge\tau_M^n}\Phi^{\prime}(\Vert u_n^M(s)\Vert^2) \Vert u_n^M(s)\Vert_{V}^2 ds \nonumber\\
   \leq & \Phi(\Vert u_0\Vert^2)e^{CT} .
\end{align}
Thus, we have
\begin{align}\label{190630.0137.1}
  \sup_n E\int_{0}^{T\wedge\tau_M^n} \Vert u_n^M(s)\Vert_{V}^2 ds \leq \frac{1}{c_M}\Phi(\Vert u_0\Vert^2)e^{CT}.
\end{align}
$\blacksquare$
%

\vskip 0.4cm
Using Lemma \ref{190630.0142} and by a slight modification of the proof of Lemma \ref{190630.0140}, we can prove the following result.
\begin{lemma}\label{190630.0149}
Suppose hypothesis (\hyperlink{H.3}{H.3}) holds. Then for any $T, M>0$, $\beta<\frac{1}{2}$ and $1<p<2$, there exists a constant $K_{T, M, \beta, p}$ such that
\begin{align}\label{5.15}
  \sup_{n} E\left(\Vert u_n^M\Vert_{W^{\beta, p}([0,T];V^*)}^p \mathbf{1}_{\{\tau_M^n\geq T\}} \right)\leq K_{T, M, \beta, p}.
\end{align}
\end{lemma}
Now we are ready to state the tightness of $u_n$.
\begin{lemma}
Suppose hypothesis (\hyperlink{H.3}{H.3}) holds. Then for any $T>0$ and $1<p<2$,  $\{u_n\}$ is tight in $L^p([0,T];H)$ and in $C([0,T];V^*)$.
\end{lemma}
\noindent {\bf Proof}.
 As in the proof of  Lemma \ref{lemma 4.5}, we know that $L^p([0,T];V)\cap W^{\beta, p}([0,T];V^*)$ is compactly embedded into $L^p([0,T];H)$. So for any $L>0$, the set
\begin{align*}
  K_L :=\left\{u\in L^p([0,T];H): \Vert u\Vert_{L^p([0,T];V)}+\Vert u \Vert_{W^{\beta, p}([0,T];V^*)}\leq L\right\}
\end{align*}
is a compact subset of $L^p([0,T];H)$.
By Lemma \ref{190630.0143}, for any $M>0$ we have
\[ \Phi(M)P(\tau_M^n\leq T)\leq E\left[\Phi(\Vert u_n(T\wedge \tau_M^n)\Vert^2)\right]\leq \Phi(\Vert u_0\Vert^2)e^{CT} , \]
namely,
\begin{equation}\label{5.16}
P(\tau_M^n\leq T)\leq \frac{\Phi(\Vert u_0\Vert^2)e^{CT}}{\Phi(M)} .
\end{equation}
For any given small positive constant $\varepsilon$, as $\Phi(M)\rightarrow \infty$ we can choose a sufficiently large constant $M$ such that
 \begin{equation}\label{5.17}
P(\tau_M^n\leq T)\leq \frac{\varepsilon}{2}, \quad \mbox{for all}\quad n\geq 1.
\end{equation}
With the constant $M$ chosen as above, using Lemma \ref{190630.0142} and Lemma \ref{190630.0149} we have for any $L>0$ and $n\geq 1$,
\begin{align}\label{5.18}
  & P(u_n\notin K_L)\leq P(u_n\notin K_L, \tau_M^n\geq T)+P(\tau_M^n\leq  T) \nonumber\\
  =&P(u_n^M\notin K_L, \tau_M^n\geq T)+P(\tau_M^n\leq  T)\nonumber\\
  \leq & \sup_{n}P\left(\Vert u_n^M\Vert_{L^p([0,T];V)}+\Vert u_n^M \Vert_{W^{\beta, p}([0,T];V^*)}> L , \tau_M^n\geq T \right)+\frac{\varepsilon}{2} \nonumber\\
  \leq & \frac{2^{p-1}}{L^p}\sup_n E\left[\left(\Vert u_n^M\Vert_{L^p([0,T];V)}^p+\Vert u_n^M \Vert_{W^{\beta, p}([0,T];V^*)}^p  \right) \mathbf{1}_{\{\tau_M^n\geq T\}} \right]+\frac{\varepsilon}{2} \nonumber\\
  \leq & \frac{2^{p-1}}{L^p}C_M^{\prime}+\frac{\varepsilon}{2},
\end{align}
where $C_M^{\prime}$ is a constant depending on $M$ but not on $L, n$. Now we can choose constant $L$ large enough so that
$$ P(u_n\notin K_L)\leq \varepsilon$$
for all $n\geq 1$. Since $\varepsilon$ is arbitrary, we have proved the tightness of $\{u_n\}$ in the space $L^p([0,T];H)$.

The tightness of $\{u_n\}$ in $C([0,T];V^*)$ can be proved by modifying the proof of Lemma \ref{lemma 4.5} combined with the above argument. We omit it here.  $\blacksquare$

\vskip 0.3cm
In view of the tightness of $\{u_n\}$, let $u$ be any weak limit of $\{u_n\}$ on some probability space $(\widetilde{\Omega}, \widetilde{\mathcal{F}},\widetilde{P})$.
\begin{lemma}
Suppose hypothesis (\hyperlink{H.3}{H.3}) holds. Then for any $T>0$,
\begin{align}
\label{190707.1128}  \widetilde{P} \left(\int_0^T \Vert u(t)\Vert_{V}^2 dt <\infty\right) & =1 , \\
\label{190707.1543}  \widetilde{P} \bigg(\sup_{t\leq T} \Vert u(t)\Vert  <\infty\bigg) & =1 .
\end{align}
\end{lemma}
\noindent {\bf Proof}. We only prove (\ref{190707.1128}) here, the proof of (\ref{190707.1543}) is similar.
Since $\{u_n\}$ is tight in $C([0,T]; V^*)$, without loss of generality and to simplify the notation, we can assume that for $P$-a.s.,
\begin{align}
  \Vert u_n -u \Vert_{C([0,T]; V^*)}\rightarrow 0 .
\end{align}
As in (\ref{0403.2}), we have for $P$-a.s.,
\begin{align}\label{190707.1007.1}
  \Vert P_m u_n-P_m u\Vert_{C([0,T_2];V)}\xrightarrow{n\rightarrow \infty} 0 , \quad \forall\, m\in\mathbb{N} .
\end{align}
Fatou's lemma and the above convergence yield
\begin{align}\label{190707.1007.2}
  & \int_{0}^{T} \Vert u(t)\Vert^2 dt
  \leq  \int_{0}^{T} \lim_{m\rightarrow \infty} \Vert P_m u(t)\Vert^2 dt
    \leq  \liminf_{m\rightarrow \infty}\int_{0}^{T} \Vert P_m u(t)\Vert^2 dt \nonumber\\
    = &  \liminf_{m\rightarrow \infty}\liminf_{n\rightarrow \infty}\int_{0}^{T} \Vert P_m u_n(t)\Vert^2 dt
    \leq  \liminf_{n\rightarrow \infty}\int_{0}^{T} \Vert u_n(t)\Vert^2 dt .
\end{align}
From (\ref{190707.1007.2}), it follows that for any $L>0$,
\begin{align}\label{190707.1026.1}
  P\left(\int_{0}^{T} \Vert u(t)\Vert_{V}^2 dt \geq L \right)
\leq  P\left(\liminf_{n\rightarrow\infty} \int_{0}^{T} \Vert u_n(t)\Vert_{V}^2 dt \geq L \right) .
\end{align}
Let
\[ A_n := \left\{\omega\in \Omega: \int_{0}^{T} \Vert u_n(t,\omega) \Vert_{V}^2 dt \geq L -1 \right\} . \]
The continuity of probability measures gives
\begin{align}\label{190707.1026.2}
  & P\left(\liminf_{n\rightarrow\infty} \int_{0}^{T} \Vert u_n(t)\Vert_{V}^2 dt \geq L \right) \nonumber\\
\leq & P \left(\liminf_{n\rightarrow\infty} A_n\right)
\leq  \liminf_{n\rightarrow\infty} P \left( A_n\right) \nonumber\\
= &  \liminf_{n\rightarrow\infty} P \left( \int_{0}^{T} \Vert u_n(t) \Vert_{V}^2 dt \geq L - 1\right) .
\end{align}
Combining (\ref{190707.1026.1}) and (\ref{190707.1026.2}) together and applying the stopping time argument lead to
\begin{align}
  & P\left(\int_{0}^{T} \Vert u(t)\Vert_{V}^2 dt \geq L \right) \nonumber\\
\leq & \liminf_{n\rightarrow\infty} P \left( \int_{0}^{T} \Vert u_n(t) \Vert_{V}^2 dt \geq L - 1 \right) \nonumber\\
\leq & \liminf_{n\rightarrow\infty} \left[ P \left( \int_{0}^{T} \Vert u_n(t) \Vert_{V}^2 dt \geq L - 1, \tau_M^n\geq T\right) + P\left(\tau_M^n <T\right) \right]\nonumber\\
\leq & \liminf_{n\rightarrow\infty} P \left( \int_{0}^{T\wedge\tau_M^n} \Vert u_n(t) \Vert_{V}^2 dt \geq L - 1\right) + \sup_n P\left(\tau_M^n <T\right) \nonumber\\
\leq & \frac{1}{L-1} \sup_n E\int_{0}^{T\wedge\tau_M^n}\Vert u_n(t) \Vert_{V}^2 dt + \sup_n P\left(\tau_M^n <T\right) .
\end{align}
Now by  (\ref{5.9}) and (\ref{5.16}), we first let $L$ go to infinity, then let $M$ go to infinity to obtain
\begin{align}
   P\left(\int_{0}^{T} \Vert u(t)\Vert_{V}^2 dt =\infty \right) =0 .
\end{align}
Hence (\ref{190707.1128}) is proved. $\blacksquare$

\vskip 0.5cm

After the preparations above, we are ready to state the  following theorem whose proof is similar to that of Theorem \ref{theorem 4.6}.

\begin{theorem}
Suppose hypothesis (\hyperlink{H.1}{H.1}) and (\hyperlink{H.3}{H.3}) hold. Then there exists a unique global solution to (\ref{Abstract}) for any initial value $u_0\in L^2(D)$.
\end{theorem}
\section{Appendix}
\setcounter{equation}{0}

There are two nonlinear versions of Gronwall's inequality used in this paper.
\begin{lemma}\label{A.1}
Let $a, b, Y$ be nonnegative functions on $\mathbb{R}$, and there are constants $c\geq 0$, $0\leq \alpha<1$  such that
\begin{align}\label{6.1}
  Y(t)\leq c+ \int_{t_0}^t \big(a(s)Y(s)+b(s)Y(s)^{\alpha}\big) ds , \quad \forall\, t\geq t_0.
\end{align}
Then for any $t\geq t_0$,
\begin{align}\label{6.2}
  Y(t)\leq \bigg\{ & c^{1-\alpha}\exp\left((1-\alpha)\int_{t_0}^t a(s)ds\right)  \nonumber\\
  & +(1-\alpha)\int_{t_0}^t b(s)\exp\left((1-\alpha)\int_s^t a(r)dr\right)ds\bigg\}^{\frac{1}{1-\alpha}} .
\end{align}
\end{lemma}

The above lemma can be found in \uppercase\expandafter{\romannumeral 12}.9 Theorem 1 (pp 360) of \cite{MPF}  and references therein.



\vskip 0.5cm

By a slight modification of the proof of Theorem 3.1 from \cite{W}, we have the following lemma.
\begin{lemma}\label{A.2}
Let $X, a, M, c_1, c_2$ be nonnegative functions on $\mathbb{R}_{+}$, $M$ be an increasing function and $M(0)\geq 1$, and $c_1, c_2$ be integrable functions on finite intervals. Assume that for any $t\geq 0$,
\begin{align}\label{6.3}
  X(t)+a(t)\leq M(t)+\int_0^t c_1(s)X(s)ds +\int_0^t c_2(s) X(s)\log X(s) ds ,
\end{align}
and the above integrals are finite. Then for any $t\geq 0$,
\begin{align}\label{6.4}
  X(t)+a(t)\leq M(t)^{\exp(C_2(t))}\exp\left(\exp(C_2(t))\int_0^t c_1(s)\exp(-C_2(s))ds\right) ,
\end{align}
where $C_2(t):=\int_0^t c_2(s)ds $.
\end{lemma}
\noindent {\bf Proof}. For the completeness, we give the proof here. Write $\log_{+}(r):= \log(r\vee 1)$. For any fixed $T>0$, let
\[
Y(t):= M(T)+\int_0^t c_1(s)X(s)ds +\int_0^t c_2(s) X(s)\log_{+} X(s) ds, \quad t\in[0,T] .
\]
We see that $Y$ is almost surely differentiable on $[0,T]$,
\[
X(t)+a(t)\leq Y(t), \quad \forall\, t\in [0,T] ,
\]
and $Y(t)\geq 1$.  This leads to
\begin{align}\label{6.5}
  Y^{\prime}(t)= & c_1(t)X(t)+ c_2(t) X(t)\log_{+} X(t) \nonumber\\
  \leq &  c_1(t)Y(t)+ c_2(t) Y(t)\log_{+} Y(t) \nonumber\\
  = &  c_1(t)Y(t)+ c_2(t) Y(t)\log Y(t) .
\end{align}
Thus
\begin{align}\label{6.6}
  \left(\log Y\right)^{\prime}(t)\leq c_1(t)+ c_2(t)\log Y(t) .
\end{align}
Solving this ordinary differential inequality, we get for any $t\in [0,T]$,
\begin{align}\label{6.7}
  \log Y(t) \leq \exp(C_2(t))\left[\log M(T) + \int_0^t c_1(s)\exp(-C_2(s))ds\right] .
\end{align}
Therefore, we obtain
\begin{align}\label{6.8}
   & X(T)+a(T) \leq  Y(T) \nonumber\\
\leq & M(T)^{\exp(C_2(T))}\exp\left(\exp(C_2(T))\int_0^T c_1(s)\exp(-C_2(s))ds\right) .
\end{align}
By the arbitrariness of $T$, (\ref{6.4}) is proved. $\blacksquare$


\vskip 0.5cm
\noindent{\bf Acknowledgement}.  This work is partially supported by National Natural Science Foundation of China (No. 11671372, No. 11431014, No. 11721101), the Fundamental Research Funds for the Central Universities (No. WK0010000057), Project funded by China Postdoctoral Science Foundation (No. 2019M652174).

\end{document}